\documentclass[reqno,a4paper,11pt]{amsart}
\usepackage{amsmath,amsthm,amssymb,mathrsfs}
\usepackage[
left = 2.25cm,
right = 2.25cm,
top = 2.2 cm,
bottom = 2.2 cm
]{geometry}
\usepackage{graphicx,xcolor,tikz,mathtools}
\usepackage{bm}
\usepackage{amsaddr}
\usepackage{enumitem}
\usepackage{esint} % Average integral symbol
\usepackage{hyperref}

\usepackage{todonotes}
\usepackage{cite}

\makeatletter
\define@key{todonotes}{Helmut}[]{%
	\setkeys{todonotes}{author=\textbf{Helmut},inline,color=red!10}}%
\define@key{todonotes}{Andrea}[]{%
	\setkeys{todonotes}{author=\textbf{Andrea},inline,color=green!20}}%
\define@key{todonotes}{Yadong}[]{%
	\setkeys{todonotes}{author=\textbf{Yadong},inline,color=cyan!20}}%
\makeatother

\theoremstyle{plain}
\newtheorem{theorem}{Theorem}[section]

\newtheorem{definition}[theorem]{Definition}
\newtheorem{lemma}[theorem]{Lemma}
\newtheorem{proposition}[theorem]{Proposition}
\newtheorem{remark}[theorem]{Remark}
\theoremstyle{definition}

\numberwithin{equation}{section}

\newcommand{\bn}{\mathbf{n}}

\renewcommand{\d}{\mathrm{d}}
\newcommand{\dx}{\,\d x}
\newcommand{\dt}{\,\d t}

\newcommand{\dtau}{\,\d \tau}

\newcommand{\norm}[1]{\left\Vert #1 \right \Vert}

\newcommand{\normmm}[1]{\left\vert #1 \right\vert}

\def\Div{\mathrm{div}\,}

\newcommand{\dist}{\mathrm{dist}}

\allowdisplaybreaks[4]

\newcommand{\R}{\mathbb{R}}

\renewcommand{\d}{\mathrm{d}}

\def\Div{\mathrm{div}\,}

\def\Ld{L^2(\Omega)}
\def\bLd{\mathbf L^2(\Omega)}

\def\N{\mathbb N}
\def\R{\mathbb R}
\def\vphi{\varphi}
\def\Hu{H^1(\Omega)}

\begin{document}
	
	\title[Asymptotic stabilization of weak solutions to phase-field equations]{Asymptotic stabilization of weak solutions \\ to phase-field equations   \\ with non-degenerate mobility\\ and singular potential}
    % \title[Navier--Stokes Equations on an Elastic Surface]{Local Well-Posedness of the Surface Navier--Stokes Equations on Elastic Membranes}
	
\author[M. Grasselli and A. Poiatti]{
	\small
Maurizio Grasselli$^\ast$ and
	Andrea Poiatti$^\ddagger$
}

\address{
	$^\ast$Dipartimento di Matematica, Politecnico di Milano, 20133 Milano, Italy
}
\email{maurizio.grasselli@polimi.it}

\address{
	$^\ddagger$ Dipartimento di Scienze Matematiche, Fisiche e Informatiche,\\
	Universit\`a degli Studi di Parma, 43124 Parma, Italy
}
\email{andrea.poiatti@unipr.it}

	\date{\today}
	
	\subjclass[2020]{35B36, 35B40, 35D30, 35Q35, 76T99}
	
    \keywords{Cahn--Hilliard equations, Allen--Cahn equation, non-degenerate mobility, strict separation property, \L ojasiewicz--Simon inequality, convergence to equilibrium.}

	\begin{abstract}
   A common paradigm in phase-field models with singular potentials is that global-in-time weak solutions converge to a single equilibrium only after undergoing asymptotic  regularization. However, in arXiv:2510.17296 we  introduced a novel method to establish the convergence to a single equilibrium for solutions to Cahn--Hilliard equations, and some related coupled systems, with non-degenerate mobility and singular potentials, under very general assumptions: we only require the existence of a global weak solution satisfying an energy inequality and then we make use of a \L ojasiewicz--Simon inequality.
   Here we take a non-trivial step further. We relax the assumptions needed to prove the precompactness of trajectories, which is an essential ingredient of the complete proof. Thanks to this generalization, we can handle all the main phase-field models, with fully general singular potentials, in a three-dimensional domain, whose asymptotic behavior has so far remained an open problem. Namely, we apply the new method to the Cahn--Hilliard equation with nonlinear diffusion, the conserved Allen--Cahn equation, and the nonlocal Cahn--Hilliard equation. In the case of second-order equations, De Giorgi's iteration argument is crucial to show that weak solutions  stay asymptotically uniformly away from pure phases (strict separation property), which is the key ingredient to apply the \L ojasiewicz--Simon inequality. We expect this new technique to have a wider range of applications to coupled problems, including hydrodynamic models like the conserved Allen--Cahn--Navier--Stokes system or the nonlocal Abels--Garcke--Gr\"un system with non-degenerate mobility. Further applications to multi-component models are also possible.

    \end{abstract}
	
	\maketitle
	
	% 	\setcounter{tocdepth}{1}
	% 	\tableofcontents
	
\section{Introduction}
\noindent
   In \cite{GPoi} we introduced a new method to prove the convergence to a single stationary state of weak solutions to Cahn--Hilliard type equations with singular potentials (e.g. of Flory--Huggins type). This method is based on a \L ojasiewicz--Simon inequality, but it is very general since it only requires the existence of a global weak solution that satisfies an energy inequality. Therefore, it can be applied for instance to the 3D Cahn--Hilliard equation with non-degenerate mobility as well as to the 3D Cahn--Hilliard--Navier--Stokes system
    with unmatched densities and viscosities proposed in \cite{AGG} (a.k.a. the AGG system) with non-degenerate mobility. We recall that the first result of this type for a Cahn--Hilliard equation with singular potential was obtained in \cite{AW}. In order to apply the \L ojasiewicz--Simon inequality, which requires the (singular) potential to be real analytic, the authors needed to prove that the solution stays uniformly away (i.e. it is strictly separated) from the pure phases asymptotically. This  was obtained by showing that the set of stationary states of each trajectory is uniformly separated and the solution is regular enough to be close to the $\omega$-limit set with respect to the maximum norm. Hence, the argument applies to solutions which suitably regularize in finite time. However, these regularization properties are not known to hold if, for instance, the Cahn--Hilliard equation has a non-constant mobility in dimension three. Another relevant case is the well-known Cahn--Hilliard--Navier--Stokes system. In \cite{Abels} the author establishes the convergence to a single equilibrium in the case of matched densities and constant mobility in dimension three. The argument is quite a \textit{tour de force} since first one needs to show regularity properties of the Cahn--Hilliard equation with a given drift possessing Leray regularity. Then one has to exploit the regularization of the Navier--Stokes system for large times, crucially hinging upon a weak-strong uniqueness principle. On the contrary, our method can be applied even to the AGG system with non-constant mobility in dimension three.  Let us briefly recall how it works focusing on the classical Cahn--Hilliard equation with non-constant mobility and Flory--Huggins potential. Consider a bounded domain $\Omega \subset \mathbb{R}^d$, $d\in\{1,2,3\}$, occupied by a binary alloy and denote by $\vphi$ the relative concentration of the difference of the two components. Then consider the following initial and boundary value problem
    \begin{align}
     \label{CH1} &\partial_t\vphi=\Div(m(\vphi)\nabla \mu),\quad\text{ in }\Omega\times(0,+\infty)\\&
      \label{CH2}  \mu= F^\prime(\vphi) - D(\vphi) ,\quad\text{ in }\Omega\times(0,+\infty).
    \end{align}
    Here, $m(\cdot)$ is a non-degenerate mobility, $\mu$ is the chemical potential (i.e., the first variation of the free energy functional with respect to $\varphi$), and $F$, the opposite of the Boltzmann--Gibbs mixing entropy density $S$, is given, for instance, by
\begin{equation}
\label{FHP}
    F(s)=\frac{\theta}{2}((1+s)\text{ln}(1+s)+(1-s)\text{ln}(1-s))
    \quad\forall\,s\in(-1,1).
\end{equation}
    Here $\theta>0$ stands for the absolute temperature. $D(\varphi)$ is the decoupling term which is usually defined as follows
    \begin{equation}
    D(\vphi)=\frac{\theta_0}{2}\vphi^2,
    \end{equation}
    where
    $\theta_0>0$ is related to the critical temperature below which phase separation takes place.
    Equations \eqref{A1}-\eqref{A2} are subject to the following boundary and initial conditions
 \begin{align}
        &\partial_n\vphi=0,\quad\text{ on }\partial\Omega\times(0,+\infty) \label{CH3}\\&
      \partial_n\mu=0,\quad\text{ on }\partial\Omega\times(0,+\infty)\label{CH4}\\&
        \vphi(0)=\vphi_0,\quad\text{ in }\Omega.\label{CH5}
    \end{align}
    These conditions ensure the conservation of the total mass, that is, $\int_\Omega \vphi(x,t)\dx$ is constant
     for all times $t\geq 0$.

Let us describe our strategy.

We consider a global weak solution $\vphi$ to \eqref{CH1}-\eqref{CH5} such that,
for any $t\geq 0$ and almost any $0\leq s\leq t$, with $s=0$ included, there holds
    	\begin{align}
E(\vphi(t))+\int_s^t\int_\Omega m(\vphi(\tau))\normmm{\nabla\mu(\tau)}^2\dx\d\tau\leq E(\vphi(s)).\label{energyineqA}
    	\end{align}
Then, we define a weaker notion of $\omega$-limit
\begin{align*}
\omega(\varphi)=\{\widetilde{\varphi}\in  H^1(\Omega):\exists t_n\to +\infty \text{ such that }\varphi(t_n)\rightharpoonup  \widetilde{\varphi}\ \text{ in }H^{1}(\Omega)\},
\end{align*}
which is nonempty because $\vphi\in BC([0,+\infty);H^1(\Omega))$ (see energy inequality). We can then prove that trajectories are precompact in $H^1(\Omega)$ and replace the weak convergence in the definition of $\omega(\vphi)$ with the strong convergence, showing also that $\omega(\vphi)$ contains only stationary points. However, in this case we only have
\begin{align}
\lim_{t\to +\infty}\dist_{H^1(\Omega)}(\vphi(t),\omega(\vphi))=0,
    \label{fundam1}
\end{align}
but this is not sufficient to apply the argument developed in \cite{AW}, because $H^1(\Omega)$ clearly does not embed into $C^0(\overline{\Omega})$ for $d=2,3$. At this point, we introduce a further idea, that is, to replace the asymptotic validity of the strict separation property with a convergence of the Lebesgue measure of a measurable set. More precisely, provided we have that $\omega(\vphi)$ is uniformly strictly separated from pure phases, we can use \eqref{fundam1} (but the $L^2(\Omega)$-distance is enough) to show that there exists $\delta_1>0$ such that
\begin{align}
\normmm{\{x\in\Omega:\ \normmm{\vphi(x,t)}\geq 1-\delta_1\}}\to 0,\quad\text{ as }t\to +\infty.
    \label{measure}
\end{align}
Here $\vert E\vert$ stands for the Lebesgue measure of a measurable set $E\subset \mathbb{R}^d$.

The fundamental idea is then to use a partition of the time interval $(T,+\infty)$, for $T>0$ large enough, into two measurable sets. More precisely, for $M>0$, we introduce
$$
A_M(T):=\{t\in[T,+\infty):\ \norm{\nabla\mu(t)}_{\bLd}\leq M\}.
$$
Here and in the sequel, if $E$ is a Lebesgue-measurable set in $\mathbb{R}^d$, then $\vert E \vert$ stands for its $d$-Lebesgue measure.
This is called the set of ``good" times and it has infinite measure (thanks to \eqref{energyineqA}), while its complement of finite measure $[T,+\infty)\setminus A_M(T)$ is the set of ``bad" times when the dissipative term is either not under control or, in a zero measure subset, it is not defined.

Exploiting \eqref{measure}, we show that the uniform strict separation property holds for all $t\in A_M(T)$, namely, there exists $\delta>0$ such that
\begin{align}
\sup_{t\in A_M(T)}\norm{\vphi(t)}_{L^\infty(\Omega)}\leq 1-\delta.
    \label{sepa2}
\end{align}
Finally, using the \L ojasiewicz-Simon inequality in $A_M(T)$, the energy inequality \eqref{energyineqA} and a simple yet crucial control of the energy in the set of hyper-dissipative bad times, we recover an inequality like the one in \cite{FS} (see Lemma \ref{Feireisl} in the Appendix) for $\nabla\mu$ which entails $\nabla \mu\in L^1(\widetilde T,+\infty;\bLd)$, for some $\widetilde T \geq T$ sufficiently large. This is enough to conclude that $\omega(\varphi)$ is a singleton.

The described method is very robust, and has also been recently applied to more complicated coupled problems like a Cahn--Hilliard--Darcy--Forchheimer system with surfactant \cite{CHFD}. Nevertheless, there is a  bottleneck in the proof of the precompactness of trajectories. Indeed, it is based entirely on the monotonicity in time of the energy. However, when the energy is more involved (see, for instance, \eqref{e} below), this argument may not suffice to retrieve precompactness. Moreover, there are cases in which the \L ojasiewicz-Simon inequality requires precompactness of trajectories in $H^2(\Omega)$ which cannot be ensured by a weak solution. The new idea proposed in this contribution is to exploit further the definition of good times. More precisely, we introduce a subset of $\omega(\vphi)$, called the set of ``good equilibrium points", which is composed only of equilibria for which there exists a convergent sequence $\{\vphi( t_n)\}$ with good times $\{t_n\}\subset A_M(T)$, for $T>0$, $t_n\to\infty$. Namely, we set
     \begin{align}
 	\omega_{g}(\vphi):=\{\vphi_*\in \omega(\vphi):\ \exists t_n\to +\infty \text{ s.t. }\{t_n\}\subset A_M(T)\text{ and }\vphi(t_n)\rightharpoonup \vphi_* \text{ in }H^1(\Omega)\}
 	\label{good0}
 \end{align}
 and we prove that trajectories restricted to good times are precompact in $H^2(\Omega)$, in the sense that, given a sequence $\{t_n\}\subset A_M(T)$, for some fixed $T>0$, there always exists a subsequence such that $\vphi(t_n)\to \vphi_*\in \omega_g(\vphi)$ strongly in $H^2(\Omega)$. This is possible since the trajectories restricted to good times are not only asymptotically uniformly separated from pure phases, as shown in \cite{GPoi}, but also enjoy elliptic regularization. With this new ingredient at hand we can now use the argument developed in \cite{GPoi}, and apply the \L ojasiewicz-Simon inequality only on good times, for which the large time $H^2(\Omega)$-closeness of trajectories to an equilibrium point (or, better, a \textit{good} equilibrium point) holds.

Let us focus on the following general initial-boundary value problem for a Cahn--Hilliard/Allen--Cahn equation
     \begin{align}
     \label{A1} &\partial_t\vphi=\alpha\Div(m(\vphi)\nabla \mu) -\beta (\mu-\overline{\mu}),\quad\text{ in }\Omega\times(0,+\infty)\\&
      \label{A2}  \mu=-\gamma\Div(a(\varphi)\nabla\vphi)+\gamma\frac{a'(\vphi)}{2}\normmm{\nabla \vphi}^2 + F^\prime(\vphi) - \sigma_1 D(\vphi) - \sigma_2 J\ast \vphi,\quad\text{ in }\Omega\times(0,+\infty)\\
        &\gamma a(\varphi)\nabla\vphi\cdot\bn=0,\quad\text{ on }\partial\Omega\times(0,+\infty) \label{Af0}\\&
        \alpha m(\varphi)\nabla\mu\cdot \bn=0,\quad\text{ on }\partial\Omega\times(0,+\infty)\label{Af1}\\&
        \vphi(0)=\vphi_0,\quad\text{ in }\Omega,\label{Af}
    \end{align}
    where $\alpha,\beta,\gamma, \sigma_j$, $j=1,2$, are non-negative constants, $\bn$ is the outward unit normal to $\partial\Omega$, and we have set some other constants equal to the unity.
    Note that the conservation of the total mass still holds.
It is worth observing that Cahn--Hilliard/Allen--Cahn equations similar to \eqref{A1}-\eqref{A2} with $\sigma_2=0$ were proposed in \cite{KK} as a simplification of a mesoscopic model for multiple microscopic mechanism in surface processes (see also \cite{KN} for the theoretical analysis).

If $\beta=\sigma_2=0$ and $\alpha=\sigma_1=1$, then \eqref{A1}-\eqref{A2} is the so-called Cahn--Hilliard equation with non-convex nonlinear diffusion (see \cite{SchimpernaPawlow} and references therein). The conserved Allen--Cahn equation (see, e.g., \cite{GPCAC} and its references, see also \cite{HJMNK}) can be recovered by setting $a=1$,
$\alpha=\gamma=\sigma_2=0$, and $\sigma_1=1$, while we get the nonlocal Cahn--Hilliard equation with $\beta=\gamma=\sigma_1=0$ and $\alpha=\sigma_2=a=1$.

The first case is the one which requires the $H^2(\Omega)$-precompactness of trajectories (see \cite[Theorem 1.1]{CGGS}) and here, for the first time, we are able to show that, for any general continuous mobility $m(\cdot)$, and non-convex nonlinear diffusion $a(\cdot)$, any weak solution converges to a unique equilibrium
under minimal assumptions and also in three-dimensional bounded domains.

Concerning the conserved Allen-Cahn equation, we know that weak solutions do regularize in finite time (see \cite{GPCAC}) but here we aim at giving minimal conditions allowing  us to treat more general settings, like, for instance, the multi-component case with general mobility matrices, for which regularization is still an open problem (see \cite{GPCAC}). The method
developed in \cite{GPoi} could be actually applied without substantial changes, because the precompactness in $H^1(\Omega)$ of the entire trajectories can be shown. On the other hand, the notion of good equilibrium points allows, for instance, to handle Allen--Cahn--Navier--Stokes systems (see, e.g., \cite{AMP, GGPC}). In this case, due to the lower regularity of  $\mu$, we cannot \textit{ a priori } show that the $L^2(\Omega)$-norm of the (volume averaged) velocity converges to zero as $t\to +\infty$, that is necessary to show the $\Hu$-precompactness (cf. \cite[Theorem 3.9]{GPoi}). On the contrary, the present argument is applicable, as it only depends on the definition of the chemical potential $\mu$ and the elliptic regularization of the system restricted to good times. We shall discuss this coupled problem in a forthcoming contribution (see \cite{GHP}). As a further result, we are also able to show that any weak solution to the conserved Allen--Cahn equation asymptotically separates from pure phases, for any potential $F$ with arbitrary singularities in $\pm 1$. This is in sharp contrast with what can be obtained in the Cahn--Hilliard equation (see \cite[Lemma 3.4]{GPoi}). Indeed, in that case the asymptotic strict separation for weak solutions only holds on the good times, but no information is derived about the bad times. On the other hand, the second-order nature of Allen--Cahn type equations allows to overcome this issue by means of De Giorgi's iterations. We point out that this is the first result on the validity of strict separation property which does not first require any regularization of the weak solution, as here we do not need the crucial condition $F'(\varphi)\in L^\infty(\tau,T;L^1(\Omega))$ (see \cite{GPCAC,GGPC}), for some $0\leq \tau<T$.

In the case of nonlocal Cahn--Hilliard equation with non-degenerate mobilities, we observe that regularization in finite time is an open problem so far. Thus, this improvement based on the regularization on good times and on the good equilibrium points seems to be necessary. Of course here the $\omega$-limit set \eqref{good0} is defined with respect to the $L^2(\Omega)$-topology. A first delicate point is the fact that the $\omega$-limit associated to a weak solution is only bounded in $H^1(\Omega)$, which is not embedded in $L^\infty(\Omega)$ in dimension two or three. Therefore, we cannot prove any uniform strict separation property in a standard way. Resorting again to the notion of good times we are able to overcome this issue and show that the $\omega$-limit is indeed uniformly strictly separated.
%As a consequence, if we restrict our attention to the good equilibrium points, then, thanks to the definition of good times, a given weak solution is uniformly strictly separated from pure phases.
Then, a careful application of De Giorgi's iterations allows us to show that this property is enough to prove that the entire weak solution is asymptotically strictly separated from pure phases.

This result is a complete novelty in the literature, since it holds for fully general singular potentials $F$. Indeed, the separation property based on the regularization in finite time is known to hold for weak solutions only under some further assumptions on the singular behavior of $F'$ at $\pm1$ (see \cite{P,GalP} for the weakest assumptions known so far). Here we show that \textit{any} global weak solution asymptotically separates from pure phases \textit{regardless} of the type of singularity of $F'$ at pure phases. After proving this property and observing that the trajectories restricted to good times are precompact in $L^2(\Omega)$, we can argue as in \cite{GPoi} and apply the \L ojasiewicz--Simon inequality, only along the good times, to prove that any global weak solution asymptotically converges to a unique equilibrium.

It is also worth observing that the present approach could be used to show a similar result for the nonlocal AGG system studied in \cite{Frigeri} (see also \cite{GGGP} for the only result known so far on the asymptotic behavior in dimension two, with constant mobility, still based on regularization of weak solutions).

The plan of the paper goes as follows. The next section contains  notation, functional setting, and the statements of the three known theorems
about the existence of a global weak solution satisfying the energy inequality. The first is concerned with the Cahn--Hilliard equation with non-convex nonlinear diffusion, the second is concerned with the conserved Allen--Cahn equation, and the third is about the nonlocal Cahn--Hilliard equation. Section~\ref{main} contains our main results about the convergence to a single equilibrium for the previous equations. These results are proven in the final Section~\ref{proofs}.

    \section{Preliminaries}
    \label{pre}
    \subsection{Notation and functional setting}
Here we introduce some notation together with the functional spaces used in the sequel.

 \begin{enumerate}[label=\textnormal{(N\arabic*)},leftmargin=*]
 	
 	\item \textbf{Notation for general Banach spaces.}
 	For any (real) normed space $X$%of scalar-valued functions
    , we denote its norm by $\|\cdot\|_X$,
 	its {dual space by $X'$}.
 	If $X$ is a Hilbert space, we write $(\cdot,\cdot)_X$ to indicate the corresponding inner product.
 	Moreover, the relative spaces of vector-valued or matrix-valued functions with each component in $X$ are denoted by $\mathbf{X}$.
 	
 	\item \textbf{Lebesgue and Sobolev spaces.}
 	Assume $\Omega$ to be a bounded domain in $\R^d$, $d=2,3$ of class $C^2$.
 	For $1 \leq p \leq \infty$ and $k \in \N$, the standard Lebesgue and Sobolev spaces defined on $\Omega$ are denoted by $L^p(\Omega)$ and $W^{k,p}(\Omega)$, and their canonical norms are denoted by $\|\cdot\|_{L^p(\Omega)}$ and $\|\cdot\|_{W^{k,p}(\Omega)}$, respectively.
 	In the case $p = 2$, we set $H^k(\Omega) = W^{k,2}(\Omega)$. The $L^2(\Omega)$ inner product is simply denoted by $(\cdot,\cdot)$.
 	Also, for any interval $I\subset\R$, any Banach space $X$, $1 \leq p \leq \infty$ and $k \in \N$, we use $L^p(I;X)$, $W^{k,p}(I;X)$ and $H^{k}(I;X) = W^{k,2}(I;X)$ to indicate the Lebesgue and Sobolev spaces of functions with values in $X$. The canonical norms are indicated by $\|\cdot\|_{L^p(I;X)}$, $\|\cdot\|_{W^{k,p}(I;X)}$ and $\|\cdot\|_{H^k(I;X)}$, respectively.
 	    We also define
\begin{align*}
    L^p_\mathrm{loc}(I;X)
    &:=
    \big\{
        u:I\to X \,\big\vert\, u \in L^p(J;X) \;\text{for every compact interval $J\subset I$}
    \big\}
    \\[1ex]
    L^p_\mathrm{uloc}(I;X)
    &:=
    \left\{ u:I\to X \,\middle|\,
    \begin{aligned}
    &u \in L^p_\mathrm{loc}(I;X) \;\text{and}\; \exists\, C>0\; \sup_{t\in I}\|u\|_{L^p(t,t+1;X)} \le C
    \end{aligned}
    \right\}.
\end{align*}
The spaces $W^{k,p}_\mathrm{loc}(I;X)$, $H^k_\mathrm{loc}(I;X)$, $W^{k,p}_\mathrm{uloc}(I;X)$, $H^k_\mathrm{uloc}(I;X)$ are defined in a similar way.
 %	For simplicity, we just write $(\cdot,\cdot) := (\cdot,\cdot)_{L^2(\Omega)}$, $\|\cdot\|:=\|\cdot\|_{L^2(\Omega)}$.
 	\item \textbf{Spaces of continuous functions.}
 	For any interval $I\subset\R$ and any Banach space $X$, $C(I;X)$ denotes the space of continuous functions mapping from $I$ to $X$ and $BC(I;X)$ denotes the space of bounded functions in $C(I;X)$. Furthermore, $C_\mathrm{w}(I;X)$ denotes the space of functions mapping from $I$ to $X$, which are continuous on $I$ with respect to the weak topology of $X$, and $BC_\mathrm{w}(I;X)$ indicates the space of bounded functions in $C_\mathrm{w}(I;X)$. Then, we denote by $C^\gamma(I;X)$, $\gamma\in(0,1]$, the space of $\gamma$-H\"{o}lder (Lipschitz, if $\gamma=1$) continuous functions with values in $X$. In addition, $C_0^k(I;X)$ stands for the space of $k$-continuously differentiable functions with compact support mapping $I$ into $X$.
 	\item \textbf{Spaces of functions with zero integral mean.}
If $v\in H^1(\Omega)'$, then we define its generalized spatial mean as
\begin{equation*}
    \overline{v}:=\frac{1}{\normmm{\Omega}} \langle v,1 \rangle_{\Hu',\Hu}.
\end{equation*}
If $v\in L^1(\Omega)$, then this spatial mean becomes the standard integral average.
By means of this definition, we introduce the following Hilbert spaces:
\begin{align*}
    L^2_{(0)}(\Omega) &:= \big\{ u\in L^2(\Omega) \,:\, \overline u = 0 \big\} \subset L^2(\Omega),\\
    H^1_{(0)}(\Omega) &:= \big\{ u\in H^1(\Omega) \,:\, \overline u = 0 \big\} \subset H^1(\Omega),\\
     H^{1}_{(0)}(\Omega)' &:= \big\{ u\in H^1(\Omega)' \,:\, \overline u = 0 \big\} \subset H^1(\Omega)'.
\end{align*}
In conclusion, for $k\in \R$, we define the affine space
$$
H^1_{(k)}(\Omega)=H^1_{(0)}(\Omega)+k.
$$
\end{enumerate}

% \texttt{	\item \textbf{Spaces of divergence-free functions.}
% We define the closed linear subspaces
% 	\begin{align*}
% 		\mathbf L^p_\sigma(\Omega)
% 		&:=\overline{\{\mathbf{u}\in \mathbf{C}^\infty_0(\Omega) \,\big\vert\, \operatorname{div}\ \mathbf{u}=0\}}^{\mathbf{L}^p(\Omega)}
% 		\subset \mathbf L^p(\Omega), \quad p\in[2,\infty),\\
% 		\mathbf H^1_\sigma(\Omega)&:= \mathbf L^2_\sigma(\Omega) \cap \mathbf H^1(\Omega).
% 	\end{align*}
%In both cases, Korn's inequality yields
%\begin{equation}
%\Vert \mathbf{u}\Vert \leq \sqrt{2}\Vert D\mathbf{u}%
%\Vert\leq \sqrt{2}\Vert \nabla \mathbf{u}\Vert
%\quad \text{ for all } \mathbf{u}\in \mathbf H^1_\sigma (\Omega).
%\label{korn}
%\end{equation}
%As a trivial consequence, $\|\nabla\cdot\|$ is a norm on $\mathbf H^1_\sigma(\Omega)$ that is equivalent to the standard norm $\|\cdot\|_{\mathbf{H}^1(\Omega)}$.
% \end{enumerate}}

    \subsection{Main assumptions}
    We enumerate here all the assumptions needed to establish the results of this contribution.
\begin{enumerate}[label=\textnormal{(A\arabic*)},leftmargin=*]
    % \item \label{ASS:1} Let $\Omega \subset \R^d$, $d=2,3$, be a bounded domain with $C^{3}$-boundary.
    \item\label{ASS:0} The non-degenerate mobility function $m(\cdot)$ satisfies the following conditions
\begin{align*}
    m\in C([-1,1]),\quad 0<m_*\leq m(s),\quad \forall s\in[-1,1].
\end{align*}
\item The function $a:[-1,1]\to \R$ representing the nonlinear diffusion satisfies:
\begin{align}
	a\in C^2([-1,1]),\quad 0<a_*\leq a(s),\quad \forall s\in[-1,1].
	\label{Assa}
\end{align}
    \item \label{ASS:S1}
    The potential $f:[-1,1]\to \R$ can be written as follows
    \begin{equation*}
        f(s)=F(s)-\frac{\theta_0}{2}s^2, \quad \forall s\in [-1,1]
    \end{equation*}
    with a given constant $\theta_0>0$, where $F\in C([-1,1])\cap C^{2}(-1,1)$ is such that
    \begin{equation*}
    \lim_{r\rightarrow -1}F^{\prime }(r)=-\infty ,
    \quad \lim_{r\rightarrow 1}F^{\prime }(r)=+\infty ,
    \quad F^{\prime \prime }(s)\geq {\theta},
    \quad F'(0)=0
    \end{equation*}
    for all $s\in (-1,1)$ and a prescribed constant $\theta\in(0,\theta_0)$.
    Without loss of generality, we also assume $F(0)=0$ and $F'(0)=0$, entailing that $F(s)\geq 0$ for all $s\in [-1,1]$.
    For the sake of convenience, we extend $f$ and $F$ onto $\R\setminus[-1,1]$ by defining
    $f(s):=+\infty $ and $F(s):=+\infty $ for all $s\in\R\setminus [-1,1]$.
\item\label{ASS:S4} The kernel $J \in W_{loc}^{1,1}(\R^d)$, with $J({x})=J(-{x})$ for almost any $x\in\mathbb{R}^d$.
    %%%%%
    % \item \label{ASS:S2} In addition to \ref{ASS:S1}, there exists $\beta>\frac12$ such that
    % \begin{equation}
    % \frac{1}{F^{\prime }(1-2\delta )}=O\left( \frac{1}{|\ln (\delta )|^{\beta }}%
    % \right) ,\quad\text{ }\dfrac{1}{|F^{\prime }(-1+2\delta )|}=O\left( \frac{1}{%
    % |\ln (\delta )|^{\beta }}\right) .  \label{est}
    % \end{equation}
    % as $\delta\to 0^+$.
    %%%%%
    % \item\label{ASS:S3} In addition to \ref{ASS:S1}, it holds
    % \begin{alignat}{2}
    % 	\frac{1}{F^{\prime}(1-2\delta)}&=O\left(\frac{1}{\vert\ln(\delta)\vert}\right),
    %     &\quad\frac{1}{F^{\prime\prime}(1-2\delta)}&=O(\delta),
    % 	\label{F}
    %     \\
    % 	\dfrac{1}{\vert F^{\prime}(-1+2\delta)\vert }&=O\left(\frac{1}{\vert\ln(\delta)\vert}\right),
    %     &\quad\dfrac{1}{F^{\prime\prime}(-1+2\delta)}&=O\left(\delta\right).
    % 	\label{F2}
    % \end{alignat}
    % as $\delta\to 0^+$.
    % Moreover, there exists $\gamma_{0}>0$ such that $F^{\prime \prime }$ is monotonously increasing on $(-1,-1+\tilde\delta_0]$ and on $[1-\tilde\delta _{0},1)$.
\end{enumerate}
\begin{remark}\label{REM:LOG}
Note that \eqref{FHP} can be written as
\begin{equation}
    f(s)=F_\mathrm{log}(s)-\frac{\theta_0}{2}s^2 \quad\text{for all $s\in [-1,1]$},
    \label{f:LOG}
\end{equation}
with $F_\mathrm{log}(\pm 1) = \theta\ln(2)$ and
\begin{equation}
    F_\mathrm{log}(s)=\frac{\theta}{2}((1+s)\text{ln}(1+s)+(1-s)\text{ln}(1-s))
    \quad\text{for all $s\in(-1,1)$}.
\label{F:LOG}
\end{equation}
Thus, it satisfies all assumptions \ref{ASS:S1}.
\end{remark}
%The results on problem \eqref{A1AGG}-\eqref{AfAGG} need the further assumptions:
%\begin{enumerate}[label=\textnormal{(H\arabic*)},leftmargin=*]
%        \item \label{ASS:Viscosity} The viscosity $\nu\in W^{1,\infty}(\R)$ satisfies
%    \begin{align*}
%        0 < \nu_* \leq \nu(s) \leq \nu^*\qquad \text{ for all }s\in\R,
%    \end{align*}
%    for some positive constants $\nu_*,\nu^*\in \R$.
%    \item \label{ASS:3} The density $\rho$ is such that
%    \begin{align*}
%        \rho(s):=\frac{1+s}2\rho_1+\frac{1-s}2\rho_2,
%    \end{align*}
%    for $\rho_1,\rho_2>0$, and we set
%    \begin{align*}
%        \rho_*:=\min\{\rho_1,\rho_2\}>0,\quad \rho^*:=\max\{\rho_1,\rho_2\}>0.
%    \end{align*}
%\end{enumerate}

\subsection{Existence of a global weak solution}

Here we report the theorems about the existence of a global weak solution to each of the three initial-boundary value problems we have introduced (see \eqref{A1}-\eqref{Af}). These results will be enough to establish the convergence to a single equilibrium of any global weak solution. In most of the cases this convergence could not be proven even using the sole technique introduced in  \cite{GPoi}. This fact motivates the present non-trivial refinement.

\textbf{Case $\alpha>0$, $\beta=0$, $\gamma>0$, $\sigma_1=1$, $\sigma_2=0$.}
    This case corresponds to the Cahn--Hilliard equation with non-degenerate mobility and nonlinear diffusion. More precisely, we have
	\begin{align}
		\label{A1b} &\partial_t\vphi=\alpha\Div(m(\vphi)\nabla \mu),\quad\text{ in }\Omega\times(0,+\infty),\\&
		\label{A2b}  \mu=-\gamma\Div(a(\varphi)\nabla\vphi)+\gamma\frac{a'(\vphi)}{2}\normmm{\nabla \vphi}^2+ f^\prime(\vphi),\quad\text{ in }\Omega\times(0,+\infty),
	\end{align}
	subject to the boundary and initial conditions
	\begin{align}
		&\gamma a(\vphi)\nabla\vphi\cdot\bn=0,\quad\text{ on }\partial\Omega\times(0,+\infty), \label{Af0b}\\&
		\alpha m(\vphi)\nabla\mu\cdot\bn=0,\quad\text{ on }\partial\Omega\times(0,+\infty),\label{Af1b}\\&
		\vphi(0)=\vphi_0,\quad\text{ in }\Omega.\label{Afb}
	\end{align}
The existence of a global weak solution to \eqref{A1b}-\eqref{Afb} can be found in \cite{SchimpernaPawlow} (see also \cite[Thms. 1.2 and 1.3]{CGGS}). It reads as follows
\begin{theorem}
 \label{weakCH}
  Let $\Omega\subset R^d$, $d=2,3$, be a bounded domain of class $C^2$, and let assumptions \ref{ASS:0}-\ref{ASS:S1} be satisfied. If $\vphi_0 \in H^1(\Omega)$ is such that $\normmm{\vphi_0}\leq 1$ almost everywhere in $\Omega$ and $\overline\vphi_0\in(-1,1)$, then
there exists a global weak solution $(\vphi,\mu)$ to problem \eqref{A1b}-\eqref{Afb}. This means that
  \begin{align*}
     & \vphi\in BC([0,+\infty);H^1(\Omega))\cap L^2_{uloc}([0,+\infty);H^2(\Omega)),\\&
     \vphi\in L^\infty(\Omega\times(0,+\infty)):\ \normmm{\vphi(x,t)}<1,\quad\text{ for a.a. }(x,t)\in \Omega\times(0,+\infty),\\&
     \partial_t\vphi\in L^2(0,\infty;H^{1}_{(0)}(\Omega)'),\\&
     \mu\in L^2_{uloc}([0,+\infty);H^1(\Omega)),
  \end{align*}
    and
    \begin{align*}
        &\langle \partial_t\vphi,v\rangle_{\Hu',\Hu}+\alpha(m(\vphi)\nabla \mu,\nabla v)=0,\quad \forall v\in H^1(\Omega), \text{ for a.a. }t\geq 0,\\&
        \mu=-\gamma\Div(a(\vphi)\nabla\vphi)+\gamma \frac{a'(\vphi)}2\normmm{\nabla\vphi}^2 +f'(\vphi)\quad\text{ a.e. in }\Omega\times(0,+\infty),
    \end{align*}
    together with $\partial_{\bn}\vphi=0$ almost everywhere on $\partial\Omega\times(0,+\infty)$.
    Additionally, for any $t\geq 0$ and almost any $s \in [0,t]$, including $s=0$, it holds
    	\begin{align}
    		E(\vphi(t))+\alpha\int_s^t\int_\Omega m(\vphi(\tau))\normmm{\nabla\mu(\tau)}^2\dx\d\tau\leq E(\vphi(s)),\label{energyineq1}
    	\end{align}
        where
       \begin{align}
           E(v)=\frac\gamma2\int_\Omega a(v) \normmm{\nabla v}^2\dx+\int_\Omega f(v)\dx,
           \label{e}
       \end{align}
       for any $v\in H^1(\Omega)$ such that $\normmm{v}\leq 1$ almost everywhere in $\Omega$.
    \end{theorem}

% \texttt{   \begin{remark}
%     To be precise the energy inequality \eqref{energyineq} can be proven to be an identity (see, for instance, \cite[Lemma 2.4]{Schimperna}). Nevertheless, in this contribution we aim at stressing the fact that our proof requires nothing more than an energy inequality.  \label{schimp}
%    \end{remark}

\textbf{Case $\alpha=0$, $\beta>0$, $\gamma>0$, $\sigma_1=1$, $\sigma_2=0$.}
This is the case of the conserved Allen--Cahn equation and the problem reads
\begin{align}
	\label{A1b1} &\partial_t\vphi+\beta(\mu-\overline \mu)=0,\quad\text{ in }\Omega\times(0,+\infty),\\&
	\label{A2b1}  \mu=-\gamma\Delta\vphi+ f^\prime(\vphi),\quad\text{ in }\Omega\times(0,+\infty),\\
    &\gamma \nabla\vphi\cdot\bn=0,\quad\text{ on }\partial\Omega\times(0,+\infty), \label{Af0b1}\\&
	\vphi(0)=\vphi_0,\quad\text{ in }\Omega.\label{Afb1}
\end{align}

The existence of a global weak solution to \eqref{A1b1}-\eqref{Afb1} is given by (see, for instance,
\cite{GPCAC}, see also \cite{GGPC})
\begin{theorem}
\label{weakAC}
Let the assumptions of Theorem \ref{weakCH} be satisfied. Then
there exists a global weak solution $(\vphi,\mu)$ to problem \eqref{A1b1}-\eqref{Afb1}. This means that
\begin{align*}
	& \vphi\in BC([0,+\infty);H^1(\Omega))\cap L^2_{uloc}([0,+\infty);H^2(\Omega)),\\&
	\vphi\in L^\infty(\Omega\times(0,+\infty)):\ \normmm{\vphi(x,t)}<1,\quad\text{ for a.a. }(x,t)\in \Omega\times(0,+\infty),\\&
	\partial_t\vphi\in L^2(0,+\infty;L^2(\Omega)),\\&
	\mu\in L^2_{uloc}([0,+\infty);L^2(\Omega)),
\end{align*}
and
\begin{align*}
	&\partial_t\vphi+\beta (\mu-\overline\mu)=0,\quad\text{ a.e. in }\Omega\times(0,+\infty),\\&
	\mu=-\gamma\Delta\vphi +f'(\vphi),\ \ \ \text{ a.e. in }\Omega\times(0,+\infty),
\end{align*}
together with $\partial_{\bn}\vphi=0$ almost everywhere on $\partial\Omega\times(0,+\infty)$.
Additionally, for any $t\geq 0$ and almost any $s \in [0,t]$, including $s=0$, it holds
\begin{align}
	E_1(\vphi(t))+\beta\int_s^t\int_\Omega \normmm{\mu(\tau)-\overline{\mu(\tau)}}^2\dx\d\tau\leq E_1(\vphi(s)),\label{energyineq2}
\end{align}
where
\begin{align}
	E_1(v)=\frac\gamma2\int_\Omega  \normmm{\nabla v}^2\dx+\int_\Omega f(v)\dx,   \label{e1}
\end{align}
for any $v\in H^1(\Omega)$ such that $\normmm{v}\leq 1$ almost everywhere in $\Omega$.
\end{theorem}

\textbf{Case $\alpha>0$, $\beta=0$, $\gamma=0$, $\sigma_1=0$, $\sigma_2=1$.} Here we consider the nonlocal Cahn--Hilliard equation with non-degenerate mobility, namely,
\begin{align}
	\label{A1b2} &\partial_t\vphi=\alpha\Div(m(\vphi)\nabla \mu),\quad\text{ in }\Omega\times(0,+\infty),\\&
	\label{A2b2}  \mu=- J\ast \vphi +F^\prime(\vphi),\quad\text{ in }\Omega\times(0,+\infty),
\end{align}
endowed with the boundary and initial conditions
\begin{align}
	&
	\alpha m(\vphi)\nabla\mu\cdot\bn=0,\quad\text{ on }\partial\Omega\times(0,+\infty),\label{Af1b2}\\&
	\vphi(0)=\vphi_0,\quad\text{ in }\Omega.\label{Afb2}
\end{align}

The main result on the existence of a global weak solution to \eqref{A1b2}-\eqref{Afb2} can be found, for instance, in \cite[Theorems 3.4, 4.1, Proposition 4.2]{GGG}  (see also \cite[
Proposition 3.1]{DGG}) in the case $m=1$. However, the the proof can be easily extended to
the case with non-degenerate mobility as long as the mobility is continuous (see, for instance, \cite{Frigeri}). In this case, We have
\begin{theorem}
	 \label{weaknonloc}
	Let $\Omega\subset R^d$, $d=2,3$, be a bounded domain of class $C^2$. Suppose that assumptions \ref{ASS:0},\ref{ASS:S1}, and \ref{ASS:S4} are satisfied. If $\vphi_0 \in L^\infty(\Omega)$ is such that $\normmm{\vphi_0}\leq 1$ almost everywhere in $\Omega$ and $\overline\vphi_0\in(-1,1)$, then
	there exists a global weak solution $(\vphi,\mu)$ to problem \eqref{A1b2}-\eqref{Afb2}. This means that
	\begin{align*}
		& \vphi\in BC([0,+\infty);L^2(\Omega))\cap L^2_{uloc}([0,+\infty);H^1(\Omega)),\\&
		\vphi\in L^\infty(\Omega\times(0,+\infty)):\ \normmm{\vphi(x,t)}<1,\quad\text{ for a.a. }(x,t)\in \Omega\times(0,+\infty),\\&
		\partial_t\vphi\in L^2(0,\infty;H^{1}_{(0)}(\Omega)'),\\&
		\mu\in L^2_{uloc}([0,+\infty);H^1(\Omega)),
	\end{align*}
	and
	\begin{align}
		\label{weakformul}&\langle \partial_t\vphi,v\rangle_{\Hu',\Hu}+\alpha(m(\vphi)\nabla \mu,\nabla v)=0,\quad \forall v\in H^1(\Omega), \text{ for a.a. }t\geq 0,\\&
		\mu=-J\ast \vphi+F'(\vphi)\quad\text{ a.e. in }\Omega\times(0,+\infty).
	\end{align}
	Additionally, for any $t\geq 0$ and almost any $s \in [0,t]$, including $s=0$, it holds
	\begin{align}
		E_2(\vphi(t))+\alpha\int_s^t\int_\Omega m(\vphi(\tau))\normmm{\nabla\mu(\tau)}^2\dx\d\tau\leq E_2(\vphi(s)),\label{energyineq3}
	\end{align}
	where
	\begin{align}
		E_2(v)=\frac14\int_\Omega\int_\Omega J(x-y) \normmm{\vphi(x)-\vphi(y)}^2\dx\d y+\int_\Omega F(v)\dx,
           \label{e2}
	\end{align}
	for any $v\in L^\infty(\Omega)$ such that $\normmm{v}\leq 1$ almost everywhere in $\Omega$.
\end{theorem}
\begin{remark}
    	Note that the conservation of mass holds for all the problems above, that is, $\overline{\vphi}(t)=\overline{\varphi_0}$ for all $t\geq 0$.
    \end{remark}

\section{Main results}
\label{main}
  \subsection{The Cahn--Hilliard equation with non-degenerate mobility and nonlinear diffusion}\label{mainCH}
 Convergence to a unique equilibrium has been successfully addressed in \cite{CGGS} when $\Omega\subset \R^2$ and $F$ is the logarithmic potential, by exploiting the nowadays standard technique of \cite{AW} based on the (instantaneous) regularization of weak solutions. They also show that, in three dimensions, if the initial data are sufficiently smooth and close to a local energy minimizer, there exists a unique global strong solution converging to an
 equilibrium point of the free energy. As already explained in the introduction, here we aim at exploiting and extending the novel approach first introduced in \cite{GPoi}, which allows to show that any \textit{global weak solution} converges to a single equilibrium without using any regularization in finite time. This shows that this convergence holds in dimension three, for any general singular potential, and for mobilities which are only continuous (cf. \cite{CGGS} where a $C^2$ regularity is required).
 Let us consider the set of admissible initial data:
 \begin{align}
 	\mathcal{H}_k:=\left\{\varphi\in H^1(\Omega): \Vert\varphi\Vert_{L^\infty(\Omega)}\leq 1,\quad \vert\overline{\varphi}\vert= k \right\}, \label{Hk}
 \end{align}
 with $k\in[0,1)$, and fix an initial datum $\varphi_0\in \mathcal{H}_k$. Let then $\varphi$ be a global-in-time weak solution departing from $\varphi_0$, which might not be unique, whose existence is ensured by Theorem \ref{weakCH}. Following \cite{GPoi}, we introduce the weak $\omega$-limit set associated to $\varphi$, i.e.,
 \begin{align*}
 	\omega(\varphi)=\{\widetilde{\varphi}\in  \mathcal{H}_k:\exists t_n\to +\infty \text{ s.t. }\varphi(t_n)\rightharpoonup \widetilde{\varphi}\ \text{ in }H^1(\Omega)\}.
 \end{align*}
 Of course $\varphi\in BC([0,\infty);H^1(\Omega))$, and thus $\omega(\vphi)$ is non-empty. We further characterize the set $\omega(\varphi)$, showing that it is composed by equilibrium points, which are defined as follows
 \begin{definition}
 	$\varphi_\infty$ is an equilibrium point to problem \eqref{A1b}-\eqref{Afb} if $\varphi_\infty\in \mathcal{H}_k\cap H^2(\Omega)$ satisfies the stationary Cahn--Hilliard equation
 	\begin{align}
 		-\gamma\Div(a(\vphi_\infty)\nabla \vphi_\infty)+\gamma\frac{a'(\vphi)}{2}\normmm{\nabla \vphi_\infty}^2 +f^\prime(\varphi_\infty)=\mu_\infty,\quad \text{ a.e. in }\Omega,
 		\label{conv1t}
 	\end{align}
 	together with $\partial_{\bn}\varphi_\infty=0$ almost everywhere on $\partial\Omega$, for some $\mu_\infty\in \R$.
 \end{definition}
 \begin{remark}
 	Note that, given $\mu_\infty\in \R$, solutions to \eqref{conv1t} do exist (see, e.g., \cite[Section 3]{CGGS}, see also \cite{AW, GGPS}), but they might be non-unique.
 \end{remark}
 If we introduce the set of all the stationary points
 $$
 \mathcal{S}:=\left\{\varphi_\infty\in \mathcal{H}_k\cap H^2(\Omega): \varphi_\infty\text{ satisfies }\eqref{conv1t}\right\},
 $$
 then we can prove that $\omega(\vphi)\subset \mathcal{S}$. Namely, we have (see Section \ref{secconvaaa} for its proof)
 \begin{lemma}
 	\label{convaaa}
 	Let the assumptions of Theorem \ref{weakCH} hold, and assume further that $\partial\Omega$ is of class $C^3$. Then, we have
 	$$
 	\omega(\vphi)\subset \mathcal{S}.
 	$$
    Moreover, there exists $E_\infty\in \R$ such that
 	\begin{align}
 		E(\vphi_\infty)=E_\infty,\quad \forall \vphi_\infty\in \omega(\vphi).\label{E1}
 	\end{align}
 	Additionally, $\omega(\vphi)$ is bounded in $ H^3(\Omega)$, and there exists $\delta_0>0$ such that
 	\begin{align}
 		\| \vphi_\infty\|_{L^\infty(\Omega)}\leq 1-2\delta_0,\quad \forall \: \vphi_\infty\in \omega(\vphi).
 		\label{sepaglobal}
 	\end{align}
 	In conclusion, the trajectories of $\vphi(\cdot)$ are precompact in $H^s(\Omega)$, $s\in(0,1)$, and it holds the characterization
 	\begin{align}
 		\omega(\varphi)=\{\widetilde{\varphi}\in  \mathcal{H}_k:\exists t_n\to +\infty \text{ such that }\varphi(t_n)\to \widetilde{\varphi}\text{ in }H^1(\Omega)\}.\label{omegal}
 	\end{align}
    Also, $\omega(\vphi)$ is compact in $H^1(\Omega)$, as well as it holds
 	\begin{align}
 		\lim_{t\to +\infty}\dist_{ H^s(\Omega)}(\varphi(t),\omega(\vphi))=0,\quad \forall s\in (0,1).\label{convergenceA}
 	\end{align}
 \end{lemma}
 \begin{remark}
	Differently from \cite{GPoi}, due to the nonlinear diffusion, we are not able to show that \eqref{convergenceA} holds with $s=1$. \label{weakern}
\end{remark}
 We can now introduce the set of good times as in \cite{GPoi}: for a fixed $M>0$ and $T>0$, we define
 $$
 A_M(T):=\{t\geq T:\ \norm{\nabla\mu(t)}_{\bLd}\leq  M\}.
 $$
 Recalling that (see \eqref{energyineq1}) $\nabla\mu\in L^2(0,+\infty;\bLd)$, also because $m(\cdot)\geq m_*$, we deduce that $A_M(T)$ is measurable and we have
 \begin{align}
     \normmm{[0,+\infty)\setminus A_M(T)}\leq \frac 1{M^2}\int_{[0,+\infty)\setminus A_M(T) }\norm{\nabla \mu(t)}_{\bLd}^2\dt\leq \frac{E(\varphi_0)}{m_*M^2}.\label{basicest}
 \end{align}
 This entails that $\normmm{A_M(T)}=+\infty$ for any $M>0$ and any $T>0$, and thus there exists at least one sequence of times $\{t_n\}\subset A_M(T)$ such that $t_n\to +\infty$.

 The fundamental novelty of this contribution is that we also define a set of \textquotedblleft good equilibrium points" as follows:
 \begin{align}
 	\omega_{g}(\vphi):=\{\vphi_*\in \omega(\vphi):\ \exists t_n\to +\infty \text{ s.t. }\{t_n\}\subset A_M(T)\text{ and }\vphi(t_n)\rightharpoonup \vphi_* \text{ in }H^1(\Omega)\}.
 	\label{good}
 \end{align}
 Observing that $\omega_g(\vphi)\subset \omega(\vphi)\subset \mathcal S$, this set collects all the equilibria for which there exists a converging subsequence evaluated only at good times. This set is of course always nonempty thanks to \eqref{basicest} and $\vphi\in BC([0,+\infty);H^1(\Omega))$, as it is enough to choose a sequence of good times $\{t_n\}\subset A_M(T)$ for some $T>0$ and obtain a subsequence of good times such that $\vphi(t_n)\rightharpoonup \vphi_*$ weakly in $H^1(\Omega)$, giving $\vphi_*\in \omega_g(\vphi)$.  Note also that, by definition, $\omega_g(\vphi)$ does not depend on the choice of $T>0$, so that we do not need to specify this dependence. The set of good equilibria $\omega_g(\vphi)$ results to be crucial, as it enjoys much stronger properties than what we have discussed for $\omega(\vphi)$ (cf., in particular, Remark \ref{weakern}). In particular, we prove the following properties (see Section \ref{sec_prooftwoparts}).
 \begin{lemma}\label{twoparts}
 	Let the assumptions of Theorem \ref{weakCH} hold and assume further that $\partial\Omega$ is of class $C^3$.
 	%    Given the set
 	%    \begin{align}
 		%    	A_\delta(t) &:= \{x\in\Omega:\; |\varphi(x,t)|\geq 1-{\delta_1}\},\quad t\geq0 ,\label{Adelta}
 		%    \end{align}
 	%    it holds
 	%    \begin{align}
 		%    	\lim_{t\to\infty }\normmm{A_\delta(t)}\to 0,
 		%    	\label{conerg}
 		%    \end{align}
 	%    where $\delta_1>0$ is given in \eqref{sepaglobal}.
 	Then, for any $M>0$ there exists $\delta\in(0,\delta_0)$ and $T_S>0$ such that
 	\begin{align}
 		\sup_{t\in A_M(T_S)}\norm{\varphi(t)}_{L^\infty(\Omega)}\leq 1-\delta. \label{asympt}
 	\end{align}
 	Furthermore, for any $r\in(0,3)$ and any $\varepsilon>0$ there exists $T=T(\varepsilon)\geq T_S$ such that
 	\begin{align}
 		\dist_{H^r(\Omega)}(\vphi(t),\omega_g(\vphi))<\varepsilon,\quad \forall t\in A_M(T).
 		\label{precomp}
 	\end{align}
   In conclusion, the set $\omega_g(\vphi)$ is compact in $H^r(\Omega)$ for any $r\in(0,3)$.
 \end{lemma}
 \begin{remark}
 	Notice that property \eqref{precomp} means that trajectories of good times are precompact in $H^2(\Omega)$, which is the basic property needed to apply the \L ojasiewicz-Simon inequality. The proof is based on some elliptic regularization properties which can be proven in the set of good times.
 \end{remark}
 \begin{remark}
     We point out that the $C^3$ regularity of $\partial\Omega$ is here required to allow for the application of elliptic regularity results giving $H^3(\Omega)$-regularity, which gives the fundamental precompactness result \eqref{precomp}.
 \end{remark}
 As a consequence of this fundamental lemma, in Section \ref{sec:proofLoja} we can eventually prove that the $\omega$-limit is formed by a unique element, if the potential is analytic in $(-1,1)$. In particular, the following result holds
 \begin{theorem}\label{uniqueeq}
 	Let the assumptions of Theorem \ref{weakCH} hold and suppose additionally that $F$ and $a$ are real analytic in $(-1,1)$. Then any global weak solution $\vphi$ given by Theorem \ref{weakCH}, departing from the initial datum $\vphi_0\in
 	\mathcal{H}_k$, converges to a single equilibrium point $\vphi_\infty\in \mathcal S$, i.e., $\omega(\vphi)=\{\vphi_\infty\}$. In particular, it holds
 	\begin{align} \lim_{t\to + \infty}\Vert\vphi(t)-\vphi_\infty\Vert_{H^s(\Omega)}=0,
 		\label{equil}
 	\end{align}
 	for any $s\in(0,1)$.
 \end{theorem}
 \begin{remark}
 	Recall that the convergence \eqref{equil} is much stronger if we restrict ourselves to the subset of good times, as indicated by \eqref{precomp}.
 \end{remark}

  \subsection{The conserved Allen--Cahn equation}\label{mainAC}
 Here we study the longtime behavior of each (weak) trajectory of the conserved Allen--Cahn equation. Convergence to a unique equilibrium has been studied in \cite{GPCAC} exploiting the regularization of weak solutions in finite time. Here we propose an application of the arguments developed in \cite{GPoi}, which does not require the instantaneous regularization. This new methodology turns out to be essential, for instance, in the case of couplings with the Navier--Stokes system in dimension three. Additionally, here we show for the first time that weak solutions to the conserved Allen--Cahn equation asymptotically separate from pure phases, without using regularization properties. The asymptotic strict separation was first shown in two and three dimensions for fully general potentials in \cite{GPCAC} by requiring the solution regularization. Note that this property is much stronger than \eqref{asympt}, since it is not only restricted to good times, but it also holds on bad times, meaning that the Allen--Cahn equation, a second-order parabolic PDE, enjoys stronger properties than the fourth order Cahn--Hilliard counterpart.

 Let us fix again an initial datum $\varphi_0\in \mathcal{H}_k$. Then, let $\varphi$ be a global weak solution departing from $\varphi_0$, whose existence is ensured by Theorem \ref{weakAC} (actually in this case the solution is even unique, though uniqueness is not needed here). Define the $\omega$-limit as follows
 \begin{align*}
 	\omega_1(\varphi)=\{\widetilde{\varphi}\in  \mathcal{H}_k:\exists t_n\to +\infty \text{ s.t. }\varphi(t_n)\rightharpoonup \widetilde{\varphi}\ \text{ in }H^1(\Omega)\}.
 \end{align*}
Also in this case, $\varphi$ is uniformly bounded in $H^1(\Omega)$, so that $\omega_1(\vphi)$ is non-empty. We thus characterize the set $\omega_1(\varphi)$ by showing that it contains equilibrium points, which are defined by
 \begin{definition}
 	$\varphi_\infty$ is an equilibrium point to problem \eqref{A1b1}-\eqref{Afb1} if $\varphi_\infty\in \mathcal{H}_k\cap H^2(\Omega)$ satisfies the stationary  equation
 	\begin{align}
 		-\gamma\Delta\vphi_\infty +f^\prime(\varphi_\infty)=\mu_\infty,\quad \text{ a.e. in }\Omega,
 		\label{conv1tb}
 	\end{align}
 	together with $\partial_{\bn}\varphi_\infty=0$ almost everywhere on $\partial\Omega$, for some $\mu_\infty\in \R$.
 \end{definition}

 If we define the set of all the stationary points as follows
 $$
 \mathcal{S}_1:=\left\{\varphi_\infty\in \mathcal{H}_k\cap H^2(\Omega): \varphi_\infty\text{ satisfies }\eqref{conv1tb}\right\},
 $$
 then we can show that $\omega(\vphi)\subset \mathcal{S}_1$. Furthermore, as pointed out above, we can prove that each weak solution asymptotically separates from pure phases.
 In particular, we will demonstrate (see Section \ref{secconvaaa1})
 \begin{lemma}
 	\label{convaaab}
 	Let the assumptions of Theorem \ref{weakAC} hold. Then, we have
 	$$
 	\omega_1(\vphi)\subset \mathcal{S}_1.
 	$$
    Also, there exists $E_{1,\infty}\in \R$ such that
 	\begin{align}
 		E_1(\vphi_\infty)=E_{1,\infty},\quad \forall \vphi_\infty\in \omega_{1}(\vphi).\label{einf}
 	\end{align}
 	Moreover, $\omega_1(\vphi)$ is bounded in $ H^2(\Omega)$, compact in $H^1(\Omega)$, and there exists $\delta_1>0$ such that
 	\begin{align}
 		\| \vphi_\infty\|_{L^\infty(\Omega)}\leq 1-2\delta_1,\quad \forall \: \vphi_\infty\in \omega_1(\vphi).
 		\label{sepaglobal1a}
 	\end{align}
 	Also, the trajectories of $\vphi(\cdot)$ are precompact in $H^s(\Omega)$, $s\in(0,1)$, and
 	\begin{align}
 		\lim_{t\to +\infty}\dist_{ H^s(\Omega)}(\varphi(t),\omega(\vphi))=0,\quad \forall s\in (0,1).\label{convergence1a}
 	\end{align}
 	Additionally, any weak solution asymptotically separates from pure phases, i.e., there exist $\delta_{A}\in(0,\delta_1)$ and $T_1>0$ such that
 	\begin{align}
 		\sup_{t\geq T_1}\norm{\vphi(t)}_{L^\infty(\Omega)}\leq 1-\delta_A.\label{sepac}
 	\end{align}
 \end{lemma}
 \begin{remark}
 	Observe that the energy $E_1$ is the standard Cahn--Hiliard energy. Then, one can also argue as in \cite[Lemma 3.3, eq. (3.5)]{GPoi} to show that actually the precompactness of trajectories is even in $H^1(\Omega)$, i.e., \eqref{convergence1a} holds with $s=1$. Here we propose a different argument to retrieve precompactness for $s\in(0,1)$, since it is more robust in case the energy has a different structure (like, e.g., in Allen--Cahn--Navier--Stokes and Cahn--Hilliard--Navier--Stokes systems, where also kinetic energy comes into play). \label{stronger}
 \end{remark}
Let us introduce the set of good times: for a fixed $M>0$ and $T>0$, we define
 $$
 A_{M,1}(T):=\left\{t\geq T:\ \norm{\mu(t)-\overline{\mu(t)}}_{\Ld}\leq  M\right\}.
 $$
 Observe that \eqref{energyineq2} entails that $\mu-\overline\mu\in L^2(0,+\infty;\Ld)$. Thus, $A_{M,1}(T)$ is measurable and such that
 \begin{align}
     \normmm{[0,+\infty)\setminus A_{M,1}(T)}\leq \frac 1{M^2}\int_{[0,+\infty)\setminus A_{M,1}(T)}\norm{\mu(t)-\overline{\mu(t)}}_{\Ld}^2\dt\leq \frac{E_1(\varphi_0)}{M^2},
     \label{basicest1}
 \end{align}
 so that $\normmm{A_{M,1}(T)}=+\infty$ for any $M>0$ and any $T>0$.
 Then, we introduce the \textquotedblleft good equilibrium points\textquotedblright :
 \begin{align}
 	\omega_{g,1}(\vphi):=\{\vphi_*\in \omega_1(\vphi):\ \exists t_n\to +\infty \text{ s.t.} \{t_n\}\subset A_{M,1}(T)\text{ and }\vphi(t_n)\rightharpoonup \vphi_*,\text{ weakly in }H^1(\Omega)\}.
 	\label{goodt}
 \end{align}
 Note that we have $\omega_{g,1}(\vphi)\subset \omega_1(\vphi)\subset \mathcal S_1$, and that $\omega_{g,1}(\vphi)$ is nonempty (see the argument used for for $\omega_g(\vphi)$) and does not depend on  $T>0$.

 The following lemma holds (see Section \ref{sec_prooftwoparts1} for its proof)
 \begin{lemma}\label{twoparts1}
 	Let the assumptions of Theorem \ref{weakAC} hold.
 	%    Given the set
 	%    \begin{align}
 		%    	A_\delta(t) &:= \{x\in\Omega:\; |\varphi(x,t)|\geq 1-{\delta_1}\},\quad t\geq0 ,\label{Adelta}
 		%    \end{align}
 	%    it holds
 	%    \begin{align}
 		%    	\lim_{t\to\infty }\normmm{A_\delta(t)}\to 0,
 		%    	\label{conerg}
 		%    \end{align}
 	%    where $\delta_1>0$ is given in \eqref{sepaglobal}.
% 	Then, for any $M>0$ there exists $\delta\in(0,\delta_1)$ and $T_S>0$ such that
% 	\begin{align}
% 		\sup_{t\in A_M(T_*)}\norm{\varphi(t)}_{L^\infty(\Omega)}\leq 1-\delta. \label{asympt}
% 	\end{align}
 	Then, for $r\in(0,2)$, for any $\varepsilon>0$ there exists $T=T(\varepsilon)\geq 0$ such that
 	\begin{align}
 		\dist_{H^r(\Omega)}(\vphi(t),\omega_{g,1}(\vphi))<\varepsilon,\quad \forall t\in A_{M,1}(T).
 		\label{precomp1a}
 	\end{align}
    Also, the set $\omega_{g,1}(\vphi)$ is compact in $H^r(\Omega)$, for $r\in(0,2)$.
 \end{lemma}
 \begin{remark}
 As observed in Remark \ref{stronger}, arguing as in \cite{GPoi} one can also show that property \eqref{precomp1a} actually holds for any $\vphi_\infty \in \omega_1(\vphi)$ and for $t\to +\infty$, not only along the good equilibrium points.
 \end{remark}
 Therefore, in Section \ref{sec:proofLoja1} we can eventually prove that the $\omega$-limit is a singleton, if the potential is analytic in $(-1,1)$. Namely, we have the following
 \begin{theorem}\label{uniqueeq1a}
 	Let the assumptions of Theorem \ref{weakCH}, and suppose $F$ to be real analytic in $(-1,1)$. Then any global weak solution $\vphi$ given by Theorem \ref{weakCH}, departing from the initial datum $\vphi_0\in
 	\mathcal{H}_k$, converges to a single equilibrium point $\vphi_\infty\in \mathcal S_1$, i.e., $\omega_1(\vphi)=\{\vphi_\infty\}$. In particular, it holds
 	\begin{align} \lim_{t\to +\infty}\Vert\vphi(t)-\vphi_\infty\Vert_{H^s(\Omega)}=0,
 		\label{equil1a}
 	\end{align}
 	for any $s\in(0,1)$.
 \end{theorem}
 \begin{remark}
 	The convergence \eqref{equil1a} is much stronger if we restrict ourselves to the subset of good times, as indicated by \eqref{precomp1a}. Also, recalling Remark \ref{stronger}, one can also reach $s=1$ in \eqref{equil1a}.
 \end{remark}

\begin{remark}
On account of Theorems \ref{uniqueeq} and \ref{uniqueeq1a}, one could also obtain a similar result for problem \eqref{A1}-\eqref{Af} with  $\alpha,\beta,\gamma,\sigma_1$ positive constants. In particular, for the Cahn--Hilliard/Allen--Cahn equation with singular potential (cf. \cite{KK,KN}). Note that, if $\beta (\mu-\overline{\mu})$ is replaced by $\beta\mu$ like in \cite{KK,KN}, then the mass is no longer conserved, but the energy inequality gives a straightforward control of the $H^1(\Omega)$-norm of $\mu(t)$.
\end{remark}

  \subsection{The nonlocal Cahn--Hilliard equation with non-degenerate mobility}\label{mainCHnonloc}
 In the case of constant mobility, convergence to a unique equilibrium was shown in dimension two, for instance, in \cite{DGG}, whereas in three dimensions the proof has been first given in \cite{P}, based on the instantaneous strict separation property. Note that all the existing proofs are based on the instantaneous regularization of weak solutions. Here, we show that also in the case of non-degenerate mobilities each trajectory converges to a unique equilibrium, without resorting to any regularization of the global weak solutions. We will make use of the novel approach introduced in \cite{GPoi} enhanced by the notion of ``good equilibrium points''.
 The set of admissible initial data is given by
 \begin{align}
 	\mathcal{\mathcal L}_k:=\left\{\varphi\in L^\infty(\Omega): \Vert\varphi\Vert_{L^\infty(\Omega)}\leq 1,\quad \vert\overline{\varphi}\vert= k \right\}, \label{Lk}
 \end{align}
 with $k\in[0,1)$. Fix an initial datum $\varphi_0\in \mathcal{L}_k$ and denote by $\varphi$ a global weak solution departing from $\varphi_0$ (not necessarily unique), whose existence is given by Theorem \ref{weaknonloc}. Then, we define the weak $\omega$-limit set associated with $\varphi$ as follows
 \begin{align*}
 	\omega_2(\varphi)=\{\widetilde{\varphi}\in  \mathcal{L}_k:\exists t_n\to +\infty \text{ s.t. }\varphi(t_n)\rightharpoonup \widetilde{\varphi}, \ \text{ weakly in }L^2(\Omega)\}.
 \end{align*}
 Note that $\omega_2(\vphi)$ is non-empty since $\vphi\in BC([0,+\infty);L^2(\Omega))$. We further characterize the set $\omega_2(\varphi)$, showing that it only contains equilibrium points which are defined by
 \begin{definition}
 	$\varphi_\infty$ is an equilibrium point to problem \eqref{A1b2}-\eqref{Afb2} if $\varphi_\infty\in \mathcal{L}_k\cap H^1(\Omega)$ satisfies the stationary nonlocal Cahn--Hilliard equation
 	\begin{align}
 	-J\ast\vphi+ F^\prime(\varphi_\infty)=\mu_\infty,\quad \text{ a.e. in }\Omega,
 		\label{conv1tb1}
 	\end{align}
 	for some $\mu_\infty\in \R$.
 \end{definition}
 \begin{remark}
 Given $\mu_\infty\in \R$, a (possibly non unique) solution to \eqref{conv1tb1} exists (see, e.g., \cite{Bates}).
 \end{remark}
 We then define the set of all the stationary points
 $$
 \mathcal{S}_2:=\left\{\varphi_\infty\in \mathcal{L}_k\cap H^1(\Omega): \varphi_\infty\text{ satisfies }\eqref{conv1tb1}\right\}
 $$
 and we prove first that $\omega_2(\vphi)\subset \mathcal{S}_2$, together with some further precompactness properties of the trajectories. More precisely, we have the following lemma (see Section \ref{secconvaaa2})
 \begin{lemma}
 	\label{convaaa3}
 	Let the assumptions of Theorem \ref{weaknonloc} hold. We have
 	$$
 	\omega_2(\vphi)\subset \mathcal{S}_2,
 	$$
    and there exists $E_{2,\infty}\in \R$ such that
 	\begin{align}
 		E_2(\vphi_\infty)=E_{2,\infty},\quad \forall \vphi_\infty\in \omega_{2}(\vphi).\label{E2}
 	\end{align}
 	In addition, the trajectories of $\vphi(\cdot)$ are precompact in $(H^s(\Omega))'$, $s\in(0,1]$, and
 	\begin{align}
 		\omega_2(\varphi)=\{\widetilde{\varphi}\in  \mathcal{L}_k:\exists t_n\to +\infty \text{ s.t. }\varphi(t_n)\to \widetilde{\varphi}\text{ in }L^2(\Omega)\}.\label{omegal3}
 	\end{align}
 	In addition, $\omega_2(\vphi)$ is bounded in $ H^1(\Omega)$ and thus compact in $L^2(\Omega)$, as well as it holds
 	\begin{align}
 		\lim_{t\to +\infty}\dist_{ (H^s(\Omega))'}(\varphi(t),\omega(\vphi))=0,\quad \forall s\in (0,1].\label{convergenceAA}
 	\end{align}
    \end{lemma}
    \begin{remark}
        Note that, in the present case, we are not able to prove immediately that $\omega_2(\vphi)$ is uniformly strictly separated from pure phases, since it is uniformly bounded only in $H^1(\Omega)$, apparently. A uniform bound in $L^\infty(\Omega)$ will be obtained after introducing the notion of good times (see Lemma \ref{twoparts2} below).
    \end{remark}
 Let us introduce the set of good times in this case, namely, for a fixed $M>0$ and $T>0$, we define
 $$
 A_{M,2}(T):=\{t\geq T:\ \norm{\nabla\mu(t)}_{\bLd}\leq  M\}.
 $$
 Note that, recalling \eqref{basicest}, it holds that $A_{M,2}(T)$ is measurable, and such that $\normmm{A_{M,2}(T)}=+\infty$ for any $M>0$ and any $T>0$.

 In the next lemma we will show that these sets are fundamental to prove that the $\omega$-limit set is uniformly strictly separated from pure phases. Also, we introduce again the set of \textquotedblleft good equilibrium points", since we need some further precompactness properties with respect to \eqref{convergenceAA}:
 \begin{align}
 	\omega_{g,2}(\vphi):=\{\vphi_*\in \omega_2(\vphi):\ \exists t_n\to +\infty \text{ s.t. }\{t_n\}\subset A_{M,2}(T)\text{ and }\vphi(t_n)\rightharpoonup \vphi_* \text{ in }L^2(\Omega)\}.
 	\label{good2}
 \end{align}
Then, we establish some fundamental further precompactness properties of the more regular set $\omega_{g,2}(\vphi)$ which are contained in the following lemma (see Section \ref{sec_prooftwoparts2})
 \begin{lemma}\label{twoparts2}
 	Let the assumptions of Theorem \ref{weaknonloc} hold.
 	%    Given the set
 	%    \begin{align}
 		%    	A_\delta(t) &:= \{x\in\Omega:\; |\varphi(x,t)|\geq 1-{\delta_1}\},\quad t\geq0 ,\label{Adelta}
 		%    \end{align}
 	%    it holds
 	%    \begin{align}
 		%    	\lim_{t\to\infty }\normmm{A_\delta(t)}\to 0,
 		%    	\label{conerg}
 		%    \end{align}
 	%    where $\delta_1>0$ is given in \eqref{sepaglobal}.
 	% Then, for any $M>0$ there exists $\delta\in(0,\delta_1)$ and $T_S>0$ such that
 	% \begin{align}
 	% 	\sup_{t\in A_M(T_*)}\norm{\varphi(t)}_{L^\infty(\Omega)}\leq 1-\delta. \label{asympt}
 	% \end{align}
    Then, there exists $\delta_2>0$ such that
 	\begin{align}
 		\| \vphi_\infty\|_{L^\infty(\Omega)}\leq 1-2\delta_2,\quad \forall \: \vphi_\infty\in \omega_{2}(\vphi).
 		\label{sepaglobal3}
 	\end{align}
 Furthermore, for $r\in(0,1)$, for any $\varepsilon>0$ there exists $T=T(\varepsilon)$ such that
 	\begin{align}
 		\dist_{H^r(\Omega)}(\vphi(t),\omega_{g,2}(\vphi))<\varepsilon,\quad \forall t\in A_{M,2}(T).
 		\label{precomp2}
 	\end{align}
    Also, the set $\omega_{g,2}(\vphi)$ is compact in $H^r(\Omega)$ with $r\in(0,1)$.
 \end{lemma}
 \begin{remark}
 	Property \eqref{precomp2} entails that trajectories restricted to good times are precompact in $L^2(\Omega)$, which is the basic property to apply the \L ojasiewicz-Simon inequality in the case of nonlocal energies.
 \end{remark}
By crucially exploiting the fact that $\omega_{g,2}(\vphi)$ is uniformly strictly separated, we can finally prove the asymptotic strict separation property, namely (see Section \ref{sec::proofsep} for the proof based on De Giorgi's iterations):
\begin{theorem}
Let the assumptions of Theorem \ref{weaknonloc} hold. Then any weak solution asymptotically separates from pure phases, i.e., there exist $\delta_B\in(0,\delta_2)$ and $T_2>0$ such that
 	\begin{align}
 		\sup_{t\geq T_2}\norm{\vphi(t)}_{L^\infty(\Omega)}\leq 1-\delta_B.\label{sepac1}
 	\end{align}
    \label{sep}
\end{theorem}
% \begin{remark}
%     An immediate corollary of this result is that indeed also the full $\omega$-limit $\omega_2(\vphi)$ is uniformly separated from pure phases, which is totally nontrivial to be found \textit{a priori}.
% \end{remark}
 As a consequence, in Section \ref{sec:proofLoja2} we can finally prove that the $\omega$-limit is formed by a unique element, if the potential is analytic in $(-1,1)$. More precisely, there holds
 \begin{theorem}\label{uniqueeq2}
 	Let the assumptions of Theorem \ref{weaknonloc} hold and assume further that $F$ is real analytic in $(-1,1)$. Then any global weak solution $\vphi$ given by Theorem \ref{weaknonloc}, departing from the initial datum $\vphi_0\in
 	\mathcal{L}_k$, converges to a single equilibrium $\vphi_\infty\in \mathcal S_2$, i.e., $\omega_2(\vphi)=\{\vphi_\infty\}$. In particular, it holds
 	\begin{align} \lim_{t\to +\infty}\Vert\vphi(t)-\vphi_\infty\Vert_{(H^s(\Omega))'}=0,
 		\label{equil2b}
 	\end{align}
 	for any $s\in(0,1]$.
 \end{theorem}
 \begin{remark}
 	Also in this case the convergence \eqref{equil2b} is  stronger if we restrict ourselves to the subset of good times, as it can be seen from \eqref{precomp2}.
 \end{remark}

\section{Proofs of the main results}
\label{proofs}
\subsection{\L ojasiewicz--Simon inequalities}
We first recall the classical \L ojasiewicz--Simon inequality (see \cite[Proposition 6.1]{AW}), which is valid thanks to the $C^2$-regularity of the boundary $\partial\Omega$:
	\begin{proposition}
				 Assume that $F$ satisfies \ref{ASS:S1} and is real analytic in $(-1,1)$. Let $\vphi\in  \mathcal H_k$ be such that $-1+\tilde\delta\leq \vphi(x)\leq 1-\tilde\delta$, for almost any $x\in \Omega$ {and }for some $\tilde\delta\in(0,1)$. {Furthermore, let } $\vphi_\infty\in \mathcal{S}_1$ {be fixed} such that $-1+\tilde\delta\leq \vphi_\infty(x)\leq 1-\tilde\delta$ for any $x\in\Omega$. Then, there exist $\vartheta\in \left(0,\frac{1}{2}\right]$, $\eta>0$ and a positive constant $C$ such that
				\begin{align}
					\vert{{E}_1}(\vphi)-{E}_1(\vphi_\infty)\vert^{1-\vartheta}\leq C\| \delta E_1(\varphi)\|_{\Hu'},
					\label{ener}
				\end{align}
				{provided that } $\|\vphi-\vphi_\infty\|_{\Hu}\leq \eta$,
				\label{Lojaw}
where $\delta E_1:H^1_{(k)}(\Omega)\to H^1_{(0)}(\Omega)'$ is the Frechét derivative of $E_1:H^1_{(k)}(\Omega)\to \R$ defined in \eqref{e1}.
\end{proposition}
We then report a modified \L ojasiewicz--Simon inequality that accounts for the nonlinear diffusion term in the energy $E$ defined in \eqref{e} (see \cite[Theorem 1.1]{CGGS}). In this case we need a higher-order norm to control the distance between each function and $\vphi_\infty\in \mathcal S$, namely, we have
	\begin{proposition}
				 Assume that $F$ satisfies \ref{ASS:S1} and suppose that $F$ and $a$ are real analytic in $(-1,1)$. {Let} $\vphi\in  \mathcal H_k\cap H^2(\Omega)$ {be} such that $\partial_\bn \vphi=0$ on $\partial\Omega$, and $-1+\tilde\delta\leq \vphi(x)\leq 1-\tilde\delta$, for {almost }any $x\in \Omega$ {and }for some $\tilde\delta\in(0,1)$. {Furthermore, let } $\vphi_\infty\in \mathcal{S}$ {be fixed} and such that $-1+\tilde\delta\leq \vphi_\infty(x)\leq 1-\tilde\delta$ for any $x\in\Omega$. Then there exist $\vartheta\in \left(0,\frac{1}{2}\right]$, $\eta>0$ and a positive constant $C$ such that
				\begin{align}
					&\vert{{E}}(\vphi)-{E}(\vphi_\infty)\vert^{1-\vartheta}\leq \nonumber\\
                    &C\|-\gamma\Div(a(\varphi)\nabla\vphi)+\gamma\frac{a'(\vphi)}{2}\normmm{\nabla \vphi}^2+ f^\prime(\vphi)-\overline{\gamma\frac{a'(\vphi)}{2}\normmm{\nabla \vphi}^2+f'(\vphi)}\|_{\Ld},
					\label{ener1}
				\end{align}
				{provided that } $\|\vphi-\vphi_\infty\|_{H^2(\Omega)}\leq \eta$.
				\label{Lojaw1}
\end{proposition}
In conclusion, we recall the \L ojasiewicz--Simon inequality for the nonlocal energy $E_2$ (see, e.g., \cite[Proposition 6.2]{DGG}):
\begin{proposition}
	Assume that $F$ satisfies \ref{ASS:S1} and is real analytic in $(-1,1)$. Let $\vphi\in \mathcal L_k\cap H^1(\Omega)$ be such that $-1+\tilde\delta\leq \vphi(x)\leq 1-\tilde\delta$, for almost any $x\in\overline{\Omega}$, for some $\tilde\delta\in(0,1)$ and $\vphi_\infty\in \mathcal{S}_2$. Then there exists $\vartheta\in \left(0,\frac{1}{2}\right]$, $\eta>0$ and a positive constant $C$ such that
	\begin{align}
	\vert{E}_2(\vphi)-{E}_2(\vphi_\infty)\vert^{1-\vartheta}\leq C\Vert F^\prime(\vphi)-J\ast\vphi-\overline{F^\prime(\vphi)-J\ast\vphi}\Vert_{H^1(\Omega)'},
	\label{ener2}
	\end{align}
	whenever $\Vert\vphi-\vphi_\infty\Vert_{\Ld}\leq \eta$.
	\label{Lojaw2}
\end{proposition}
\subsection{The Cahn--Hilliard equation with non-degenerate mobility and nonlinear diffusion}
 \subsubsection{Proof of Lemma \ref{convaaa}}\label{secconvaaa}
We consider a sequence $t_n\to +\infty$ such that $\varphi(t_n)\rightharpoonup  \widetilde{\varphi}$ weakly in $H^1(\Omega)$, with $\widetilde\vphi\in \omega(\varphi)$, and we can focus on a nonrelabeled subsequence such that $\vphi(t_n)\to\widetilde\vphi$ strongly in $\Ld$. We set $\varphi_n(t):=\varphi(t+t_n)$ and $\mu_n(t):=\mu(t+t_n)$, and we observe that
\begin{align}
      \label{eqs1_a}  & \langle\partial_t\vphi_n,v\rangle_{H^1(\Omega)',H^1(\Omega)}+\alpha(m(\vphi_n)\nabla \mu_n,\nabla v)=0,\quad \forall v\in H^1(\Omega), \text{ for a.a. }t\geq 0,\\&
        \mu_n=-\gamma\Div(a(\vphi_n)\nabla\vphi_n)+\gamma\frac{a'(\vphi_n)}{2}\normmm{\nabla\vphi_n(t)}^2+f'(\vphi_n)\quad\text{ a.e. in }\Omega\times(0,+\infty),\label{eqs2_a}
    \end{align}
   together with $\partial_{\bn}\varphi_n=0$ almost everywhere on $\partial\Omega \times [0,+\infty)$.
By means of Theorem \ref{weakCH}, recalling the energy inequality \eqref{energyineq1}, we get that ${{E}}(\vphi(t_n))\leq {{E}}(\vphi_0)$ for any $n\in\mathbb{N}$. Thus, for any $T>0$, there exists $C(T)>0$ independent of $n$ such that
				\begin{align}
				\| \vphi_n\|_{L^\infty(0,T;H^1(\Omega))}+\| \mu_n\|_{L^2(0,T;H^1(\Omega))}\leq C(T).\label{energy}
				\end{align}
                Note that the second term  can be controlled (see, e.g., \cite[Eq. (5.17)]{CGGS}) via energy inequality and a suitable Poincaré's inequality (see, e.g., \cite[Eq. (2.3)]{CGGS}), using the mass conservation and the control (see \cite{Kenmochi})
                $$
                \int_\Omega\normmm{F'(\vphi_n)}\dx\leq C(\overline\vphi_0)\left(1+\int_\Omega F'(\vphi_n)(\vphi_n-\overline\vphi_n)\dx\right),
                $$
                which gives (see, for the nonlinear diffusion case, \cite[Eqs. (5.17)-(5.18)]{CGGS})
                \begin{align}
                   \normmm{\int_\Omega\frac{\mu_n(t)}{\sqrt{a(\vphi_n(t))}}\dx}\leq C\int_\Omega\normmm{f'(\vphi_n)}\dx\leq C(1+\norm{\nabla\mu_n(t)}_{\bLd}).\label{ctr}
                \end{align}
Then, by comparison, we also get
\begin{align}
    \| \partial_t\vphi_n\|_{L^2(0,T;H^1(\Omega)')}\leq C(T).\label{cmp}
\end{align}
Moreover, using for instance \cite[Lemma 3.1]{CGGS} and recalling \eqref{ctr}, we obtain
$$
\norm{F'(\vphi_n)}_{L^2(0,T;L^2(\Omega))}\leq C(1+\norm{\nabla\mu}_{L^2(0,T;\bLd)})\leq C(T).
$$
Thus, by a suitable application of elliptic regularity (see \cite[Eq. (5.21)]{CGGS}), we find
\begin{align}
\| \vphi_n\|_{L^2(0,T;H^2(\Omega))}\leq C(T).
    \label{elliptic}
\end{align}
From estimates \eqref{energy}-\eqref{elliptic}, we deduce that there exists $(\vphi^*,\mu^*)$ such that, for any fixed $T>0$,
				\begin{align}
	&
        \label{est}\vphi_n\rightharpoonup \vphi^*\quad\text{ in } L^2(0,T;H^2(\Omega))\cap H^1(0,T;\Hu'),\\&
					\vphi_n\overset{\ast}{\rightharpoonup} \vphi^*\quad\text{ in } L^\infty(\Omega\times(0,T)),\label{wstar}
					\\&
					\vphi_n\to \vphi^*\quad \text{ in } L^2(0,T;H^s(\Omega)),\, \forall\: s\in[0,2)\text{ and a.e. in }\Omega\times(0,T),\label{essenziale}\\&
					\mu_n\rightharpoonup  \mu^*\quad\text{ in } L^2(0,T;\Hu).
				\end{align}
				% Additionally, by the Aubin-Lions Lemma we also infer that, for any $T>0$,
				% \begin{align}
				% 	\vphi_n\to \vphi^*\quad\text{strongly in }C([0,T];L^2(\Omega)).
				% 	\label{AL}
				% \end{align}
				These convergences are enough to pass to the limit in the equations \eqref{eqs1_a}-\eqref{eqs2_a}, meaning that {the limit pair } $(\vphi^*,\mu^*)$ satisfies, for any $T>0$,
				\begin{align*}
					\langle \partial_t\vphi^*,v\rangle_{H^1(\Omega)',H^1(\Omega)}+\alpha(m(\vphi^*)\nabla \mu^*, \nabla v)=0,\quad \forall \: v\in \Hu,\quad \text{a.e. in }(0,T),\\
					\mu^*=-\gamma\Div(a(\vphi^*)\nabla\vphi^*)+\gamma\frac{a'(\vphi^*)}2\normmm{\nabla\vphi^*}^2+f'(\vphi^*),\quad\text{a.e. in }\Omega\times(0,T),
				\end{align*}
				subject to the boundary condition $\partial_{\bn}\varphi^*=0$ almost everywhere on $\partial\Omega$ with initial datum $\vphi^*(0)={\vphi_\infty}$. The latter condition follows from the fact that $\vphi_n(0)=\vphi(t_n) \to \vphi_\infty$ strongly in $\Ld$. Concerning the convergence, we point out that the only significative difference with respect to the standard Cahn--Hilliard case (see, e.g., \cite{GPoi}) is the term multiplied by $a'$. The convergence can be shown as follows. Given $\xi\in C_c(0,\infty;C^\infty_c(\Omega))$, we have for any $T>0$
                \begin{align*}
                &\normmm{\int_0^T\int_\Omega \gamma \frac{a'(\vphi_n)}2\normmm{\nabla\vphi_n(t)}^2\xi\dx\dt-\int_0^T\int_\Omega \gamma \frac{a'(\vphi^*)}2\normmm{\nabla\vphi^*}^2\xi\dx\dt}\\&
                \leq
                \normmm{\int_0^T\int_\Omega \gamma \frac{a'(\vphi_n)}2\normmm{\nabla\vphi^*}^2\xi\dx\dt-\int_0^T\int_\Omega \gamma \frac{a'(\vphi^*)}2\normmm{\nabla\vphi^*}^2\xi\dx\dt}\\&
                \quad +\normmm{\int_0^T\int_\Omega \gamma \frac{a'(\vphi_n)}2\normmm{\nabla\vphi_n(t)}^2\xi\dx\dt-\int_0^T\int_\Omega \gamma \frac{a'(\vphi_n)}2\normmm{\nabla\vphi^*}^2\xi\dx\dt}\\&
                \leq
                \gamma \norm{a'(\vphi_n)-a'(\vphi^*)}_{L^2(0,T;L^3(\Omega))}\norm{\nabla\vphi^*}_{L^4(0,T;\mathbf L^3(\Omega))}^2\norm{\xi}_{L^\infty(\Omega\times(0,T))}\\&\quad +\gamma \norm{a'(\vphi_n)}_{L^\infty(\Omega\times(0,T))}\norm{\nabla(\vphi_*-\vphi_n)}_{L^2(\Omega\times(0,T))}\norm{\nabla(\vphi_*+\vphi_n)}_{L^2(\Omega\times(0,T))}\norm{\xi}_{L^\infty(\Omega\times(0,T)}\\&
                \quad\leq C(T)(\norm{\vphi_n-\vphi^*}_{L^2(0,T;L^3(\Omega))}+\norm{\vphi_*-\vphi_n}_{L^2(0,T;H^1(\Omega))})\norm{\xi}_{L^\infty(\Omega\times(0,T)}\to 0,
                \end{align*}
               as $n\to\infty$, where we used the regularity of $a'$, the uniform bound of $\vphi_n$ in $L^\infty(0,T;H^1(\Omega))$, \eqref{wstar}-\eqref{essenziale}, together with the embedding $ L^\infty(0,T;H^1(\Omega))\cap L^2(0,T;H^2(\Omega))\hookrightarrow L^4(0,T;L^3(\Omega))$. Furthermore, again from \eqref{wstar}-\eqref{essenziale} we deduce that
               \begin{align}
                 \label{strongc} & \vphi_n(t)\to \vphi^*(t)\quad\text{ in }H^s(\Omega),\quad\forall s\in(0,2),\\&
                    \vphi_n(t)\to \vphi^*(t)\quad\text{ in }L^p(\Omega),\quad\forall p\geq1,
               \end{align}
               as $n\to\infty$, for almost any $t\geq0$. As a consequence, by Lebesgue's dominated convergence theorem, we immediately infer
               \begin{align*}
                   \int_\Omega f(\vphi_n(t))\dx\to \int_\Omega f(\vphi(t))\dx,\quad \text{ as }n\to\infty,
               \end{align*}
               for almost any $t\geq0$. Also, we have
               \begin{align*}
                   &\normmm{\int_\Omega a(\vphi_n(t))\normmm{\nabla\vphi_n(t)}^2\dx-\int_\Omega a(\vphi^*(t))\normmm{\nabla\vphi^*(t)}^2\dx}\\&
                   \leq
                  \normmm{\int_\Omega a(\vphi_n(t))\normmm{\nabla\vphi_n(t)}^2\dx-\int_\Omega a(\vphi^*(t))\normmm{\nabla\vphi_n(t)}^2\dx}\\&\quad + \normmm{\int_\Omega a(\vphi^*(t))\normmm{\nabla\vphi_n(t)}^2\dx-\int_\Omega a(\vphi^*(t))\normmm{\nabla\vphi^*(t)}^2\dx}\\&
                  \leq
                  \norm{a(\vphi_n(t))-a(\vphi^*(t))}_{L^\infty(\Omega)}\norm{\nabla\vphi_n(t)}_{\bLd}^2\\&\quad +\norm{a(\vphi^*(t))}_{L^\infty(\Omega)}\norm{\nabla (\vphi_n(t)-\vphi^*(t))}_{\bLd}\norm{\nabla(\vphi_n(t)+{\vphi^*(t)})}_{\bLd}\\&
                  \leq C(T)\norm{\vphi_n(t)-\vphi^*(t)}_{H^\frac74(\Omega)}\to 0,
               \end{align*}
               where we used the embedding $H^\frac74(\Omega)\hookrightarrow L^\infty(\Omega)$, together with the Lipschitz regularity of $a$ and \eqref{strongc} with $s=\tfrac 74$. Therefore, we eventually infer
               \[
				\lim_{n\to \infty}{{E}}(\vphi_n(t))={{E}}(\vphi^*(t))
				\] for {almost any} $t\geq 0$. By the energy inequality \eqref{energyineq1}, we deduce that the energy ${{E}}(\vphi(\cdot))$ is nonincreasing, thus there exists ${{E}}_\infty$ such that
                \begin{align}
				\lim_{t\to+ \infty}{{E}}(\vphi(t))={{E}}_\infty.
				\label{eneA}
                \end{align}
				 Hence, {for almost any $t\geq0$}, we have
				\begin{align}
				    {{E}}(\vphi^*(t))=\lim_{n\to \infty}{E}(\vphi_n(t))=\lim_{n\to \infty}{{E}}(\vphi(t+t_n))={{E}}_\infty,\label{energyl}
				\end{align}
                so that $E(\vphi^*(\cdot))$ is constant in time and equal to $ E_\infty$ for almost any $t\geq0$. Passing then to the limit in the energy inequality, which is valid for each $\vphi_n$ thanks to \eqref{energyineq1}, we obtain
				
			\begin{align}\label{zeros}
				{{E}}_\infty+\alpha\int_s^{t}\int_\Omega m(\vphi^*(\tau))\normmm{\nabla\mu^*(\tau)}^2\dx \: \d \tau\leq {{E}}_\infty\quad \text{ for almost any } 0\leq s\leq t<\infty,\end{align}
			with $s=0$ included.

    Then, since by assumption $m(\cdot)\geq m_*>0$, \eqref{zeros} entails $\mu^*=const$ almost everywhere in $\Omega$, with a possible dependence on time. By comparison, it also holds $\partial_t\vphi^*=0$ in $\Hu'$, for almost every $t\geq 0$. As a consequence, we infer that $$\vphi^*(t)={\vphi_\infty}$$ almost everywhere in $\Omega$, for all $t\geq 0$, and thus $\mu_\infty$ is constant also in time. Therefore, ${\vphi_\infty}$ satisfies \eqref{conv1t} for some constant $\mu_\infty\in \R$, and then ${\vphi_\infty}\in \mathcal{S}$. Also, we conclude from \eqref{energyl} that
    \begin{align}
        E(\vphi_\infty)=E_\infty,\label{Einfty}
    \end{align}
and $E_\infty$ does not depend on the specific choice of $\vphi_\infty$, so that \eqref{Einfty} holds for any $\vphi_\infty\in \omega(\vphi)$, giving \eqref{einf}. Also, thanks to \eqref{strongc}, we can choose for instance $t>0$ such that $\varphi(t+t_n)\to \vphi_\infty$ strongly in $H^1(\Omega)$, entailing that the characterization \eqref{omegal} holds.

We now show the uniform strict separation properties of the $\omega$-limit, which is very similar to what obtained in \cite{GPoi}, with a slight adaptation due to the presence of the nonlinear diffusion.  It is enough to show that, given $\varphi_\infty\in \mathcal S$, there exists $\delta_{\vphi_\infty}>0$, possibly depending on $\varphi_\infty$ such that
\begin{align*}
    \norm{\vphi_\infty}_{L^\infty(\Omega)}\leq 1-\delta_{\vphi_\infty},
\end{align*}
which is trivially seen from \eqref{conv1t} (see \cite[Proposition 3.2]{CGGS}). Then, being $\omega(\vphi)\subset \mathcal S$, the same property holds for any element of the $\omega$-limit. To show that the separation property is actually uniform over  $\omega(\vphi)$, we need to show that  $\omega(\vphi)$ is bounded in $H^2(\Omega)$. First, by the energy inequality \eqref{energyineq1} we know that there is a constant $C>0$, only depending on the initial datum, such that
\begin{align}
\norm{\nabla \varphi_\infty}_{\bLd}\leq\liminf_{n\to\infty}\norm{\nabla \vphi(t_n)}_{\bLd} \leq \sup_{t\geq0}\norm{\nabla \vphi(t)}_{\bLd}\leq C(1+E(\vphi_0)),\label{f}
\end{align}
for all $\varphi_\infty\in \omega(\vphi)$.

Let us now introduce, following \cite{SchimpernaPawlow}, the quantity
$$
A(s)=\int_0^s \sqrt{a(t)}\dt.
$$
Then, we know that $\vphi_\infty\in \omega(\vphi)$ solves
\begin{align*}
   & -\gamma\sqrt{a(\vphi_\infty)}\Delta A(\vphi_\infty)+F'(\vphi_\infty)={\theta_0}\vphi_\infty+\mu_\infty,\quad \text{ a.e. in }\Omega,\\&
    \partial_\bn A(\vphi_\infty)=0\quad\text{ a.e. on }\partial\Omega.
\end{align*}
% First, since $\vphi_\infty\in[-1,1]$, from \cite[Lemma 3.1]{CGGS} we obtain
% \begin{align*}
%     \norm{F'(\vphi_\infty)}_{L^\infty(\Omega)}\leq C(\normmm{\mu_\infty}+1).
% \end{align*}
 By multiplying \eqref{conv1t} by $-\Delta A(\varphi_\infty)$, and integrating over $\Omega$, we get (recall that $\vphi_\infty$ is strictly separated from pure phases, so that the computations are rigorous)
\begin{align*}
        \int_\Omega \sqrt{a(\vphi_\infty)}\normmm{\Delta A(\vphi_\infty)}^2\dx+\int_\Omega\sqrt{a(\vphi_\infty)} F''(\vphi_\infty)\normmm{\nabla\vphi_\infty}^2\dx=\theta_0\int_\Omega  \sqrt{a(\vphi_\infty)}\normmm{\nabla \vphi_\infty}^2\dx.
\end{align*}
Observing that $\norm{a(\vphi_\infty)}_{L^\infty}\leq \max_{s\in[-1,1]}\normmm{a(s)}$ and using \eqref{f} we infer
\begin{align*}
    \norm{\Delta A(\vphi_\infty)}_{\Ld}\leq C,
\end{align*}
 which gives, since $\partial_\bn A(\vphi_\infty)=0$ almost everywhere on $\partial\Omega$,
 \begin{align}
    \norm{ A(\vphi_\infty)}_{H^2(\Omega)}\leq C. \label{ct1}
\end{align}
 Arguing as in \cite[Eq. (5.21)]{CGGS}, we thus infer
\begin{align}
    \norm{\vphi_\infty}_{H^2(\Omega)}\leq C, \label{ct2}
\end{align}
and since $C$ only depends on the initial datum, we can also write
\begin{align*}
\sup_{\vphi_\infty\in \omega(\vphi)}\norm{\vphi_\infty}_{H^2(\Omega)}\leq C,
\end{align*}
entailing that $\omega(\vphi)$ is uniformly bounded in $H^2(\Omega)$.
As a consequence, due to the compact embedding $H^2(\Omega)\hookrightarrow \hookrightarrow C^\alpha(\overline\Omega)$ for some $\alpha\in(0,1)$, by a standard contradiction argument we get the uniform strict separation property of $\omega(\vphi)$, i.e., \eqref{sepaglobal} (cf. \cite{AW} or \cite[Proof of Lemma 3.11]{GGPS}).

Thanks to the uniform separation property, we can now improve the regularity of $\omega(\vphi)$. In particular, using \eqref{f}, \eqref{ct1}, \eqref{ct2} and the uniform strict separation property \eqref{sepaglobal}, we can argue as in \cite[Proposition 3.2]{CGGS} and deduce that
\begin{align*}
    \sup_{\vphi_\infty\in \omega(\vphi)}\norm{\vphi_\infty}_{H^3(\Omega)} \leq C,
\end{align*}
i.e., that $\omega(\vphi)$ is uniformly bounded in $H^3(\Omega)$.
In conclusion, on account of \eqref{omegal}, we can write
\begin{align*}
    \omega(\vphi)=\bigcap_{t\geq 0}\overline{\bigcup_{\tau\geq t}\vphi(\tau)}^{H^1(\Omega)},
\end{align*}
so that $\omega(\vphi)$ is closed and precompact in $H^1(\Omega)$. Thus, it is compact. Finally, to show \eqref{convergenceA}, it is enough to recall that $\vphi\in BC([0,+\infty);H^1(\Omega))$, and thus the trajectories are precompact in $H^s(\Omega)$, for $s\in(0,1)$. The proof is concluded.
\subsubsection{Proof of Lemma \ref{twoparts}}\label{sec_prooftwoparts}
The proof of the validity of the asymptotic strict separation property in the set of good times can be performed in many ways, for instance by means of De Giorgi's iterations applied pointwise  to the elliptic equation that defined the chemical potential $\mu$ (see \cite{GPoi}). Here, for the sake of brevity, we propose a simpler alternative method based on elliptic regularization arguments. In particular, let us recall that, on good times, for fixed $T>0$, we have
\begin{align}
    \sup_{t\in A_M(T)}\norm{\nabla \mu(t)}_{\bLd}\leq M.\label{M1}
\end{align}
Also, by the same argument as to obtain \eqref{ctr}, since $\vphi\in[-1,1]$ and $\vphi\in L^\infty(0,+\infty;H^1(\Omega))$, we have
\begin{align*}
                   \normmm{\int_\Omega\frac{\mu(t)}{\sqrt{a(\vphi(t))}}\dx}\leq C\int_\Omega\normmm{f'(\vphi(t))}\dx\leq C(1+\norm{\nabla\mu(t)}_{\bLd}),
                \end{align*}
for almost any $t>0$. Therefore, using Poincaré's inequality \cite[Eq. (2.3)]{CGGS} and \eqref{M1}, we infer that
\begin{align}
    \sup_{t\in A_M(T)}\norm{\mu(t)}_{H^1(\Omega)}\leq C(1+M).\label{M2}
\end{align}
We can thus follow the very same elliptic argument leading to \cite[Eq. (5.21)]{CGGS} to deduce that
\begin{align*}
    \norm{\Delta\vphi(t)}_{\Ld}\leq C(1+M),\quad \forall t\in A_M(T),
\end{align*}
giving the uniform $H^2(\Omega)$-control
\begin{align}
   \sup_{t\in A_M(T)} \norm{\vphi(t)}_{H^2(\Omega)}\leq C(1+M),\label{pr}
\end{align}
for any given $T>0$.
Now it is immediate to prove that \eqref{precomp} holds for any $r\in(0,2)$. Indeed, assume that this is false. This means that there exist $r\in(0,2)$, $\varepsilon>0$, $\overline T>0$, and a sequence of good times $\{t_n \}\subset A_M(\overline T)$, $t_n\to +\infty$ as $n\to\infty$, such that
\begin{align}
    \dist_{H^r(\Omega)}(\vphi(t_n),\omega_g(\vphi))>\varepsilon,\quad \forall n\in \N. \label{contrad}
\end{align}
On the other hand, thanks to \eqref{pr}, which gives that the sequence $\{\vphi(t_n)\}$ is uniformly bounded in $H^2(\Omega)$, and the compact embedding $H^2(\Omega)\hookrightarrow\hookrightarrow H^r(\Omega)$, we deduce that there exists a subsequence $\{t_{n_k}\}_k\subset A_M(\overline{T})$ such that $\vphi(t_{n_k})\to \vphi_*$ strongly in $H^r(\Omega)$ as $k\to\infty$. Of course, by the definition of $\omega_g(\vphi)$, this entails that $\vphi_*\in \omega_g(\vphi)$. Hence, this contradicts \eqref{contrad}.
Also, this argument gives the characterization of $\omega_g(\vphi)$ as
\begin{align}
 	\omega_{g}(\vphi):=\{\vphi_*\in \omega(\vphi):\ \exists t_n\to +\infty \text{ s.t. }\{t_n\}\subset A_{M}(T)\text{ and }\vphi(t_n)\to\vphi_*,\text{ in }H^r(\Omega)\},
 	\label{good3}
 \end{align}
for any $r\in (0,2)$, and we can choose $r>\frac{d}{2}$. As a consequence, we immediately infer by \eqref{precomp}, which so far holds for $r\in(0,2)$, that $\omega_g(\vphi)$ is compact in $H^r(\Omega)$. Since $H^r(\Omega)\hookrightarrow\hookrightarrow L^\infty(\Omega)$, for $r>\tfrac d2 $, recalling the uniform separation property \eqref{sepaglobal} over $\omega(\vphi)$, there exists a finite number $N>0$ of $L^\infty$-balls $B_{\delta_0}(\vphi_m)$, centered at $\{\vphi_m\}_{m=1}^N\subset \omega_g(\vphi)$, with radius $\delta_0$, such that
$$
\omega_g(\vphi)\subset U:=\bigcup_{m=1}^NB_{\delta_0}(\vphi_m),
$$
so that it holds
\begin{align}
    \norm{v}_{L^\infty(\Omega)}\leq \norm{v-\vphi_m}_{L^\infty(\Omega)}+\norm{\vphi_m}_{L^\infty(\Omega)}\leq 1-\delta_0,\quad \forall v\in U.
\end{align}
Therefore, thanks to the fact that we have shown \eqref{precomp} for $r\in(\tfrac d2,2)$, we deduce that there exists $\tilde{T}>0$ such that $\vphi(t)\in U$ for any $t\in A_M(\tilde T)$, and thus
\begin{align}
\sup_{t\in A_M(\tilde T)}\norm{\vphi(t)}_{L^\infty(\Omega)}\leq 1-\delta_0,
    \label{unisep}
\end{align}
entailing \eqref{asympt}.

Now, thanks to \eqref{pr}, together with \eqref{asympt} and \eqref{M2}, we can argue in the same way as to obtain \cite[Eq. (7.12)]{CGGS} to conclude by elliptic regularity that
\begin{align}
   \sup_{t\in A_M(\tilde T)} \norm{\vphi(t)}_{H^3(\Omega)}\leq C(1+M).
    \label{finalcontrol}
\end{align}
Note that here the $C^3$-regularity of $\partial \Omega$ is essential.

Then, using again the same contradiction argument as above, since $H^3(\Omega)\hookrightarrow\hookrightarrow H^r(\Omega)$ for any $r\in(0,3)$, we can now extend to this range the validity of \eqref{precomp}, concluding the proof, since now the characterization \eqref{good3} holds with $r\in(0,3)$.

\subsubsection{Proof of Theorem \ref{uniqueeq}}\label{sec:proofLoja}
The proof is now similar to the one of \cite[Theorem 3.5]{GPoi}, with the fundamental difference that here we consider $\omega_g(\vphi)$ in place of the standard $\omega(\vphi)$. Let us fix $M>0$. Then we choose $\tilde \gamma$ coinciding with the value of $\delta$ given in Lemma \ref{twoparts} (see \eqref{asympt}), so that, as $\tilde \gamma\leq\delta_0$, it holds, by \eqref{sepaglobal}, $-1+\tilde\gamma\leq \vphi_\infty\leq 1-\tilde\gamma$ in $\Omega$, for any $\vphi_\infty\in\omega(\vphi)$. Furthermore, for any $\vphi_{\infty,m}\in \omega(\vphi)$ we can find $\vartheta_m\in \left(0,\frac{1}{2}\right]$ and $\eta_m>0$, given {by }Proposition \ref{Lojaw1}, for which \eqref{Lojaw1} is valid with a constant $C_m$. From Lemma \ref{twoparts} we have that $\omega_g(\vphi)\subset \omega(\vphi)$ is compact in $H^2(\Omega)$. We can then find a finite family of open $H^2(\Omega)$-balls, say $\{B_{\eta_m}\}_{m=1}^{M_1}$, centered at $\{\vphi_{\infty,m}\}_{m=1}^{M_1}\subset \omega_g(\vphi)$ and with radii $\eta_m$ (depending on the center $\vphi_{m,\infty}\in\omega_g(\vphi)$), such that
				\begin{equation*}
					\bigcup_{ \varphi_\infty \in \omega_g(\varphi)} \{\varphi_\infty\} \subset V:= \bigcup_{m=1}^{M_{1}}B_{\eta_m}.
				\end{equation*}
	Recalling \eqref{E1}, which is valid for any $\vphi_\infty\in \omega(\vphi)$, we infer that the energy functional $E(\cdot)$ is constant over $\omega(\vphi)$. Additionally, since {the centers  $\{\vphi_m\}_{m = 1}^{M_1}$ are in finite number, }we can infer that \eqref{ener1} holds \textit{uniformly}, with suitable constants, for any $\vphi\in V$ such that $\|\vphi\|_{L^\infty(\Omega)}\leq 1-\tilde\gamma$, and we can substitute $E(\vphi_\infty)$ with $E_\infty$.

    By \eqref{precomp}, there exists ${T_*}>0$ such that $\vphi(t)\in V$ for any good time $t\in A_M(T_*)$ and, additionally, by \eqref{asympt}, the uniform strict separation property holds on the set of good times, i.e.,
                \begin{align*}
                    \sup_{t\in A_M(T_*)}\norm{\vphi(t)}_{L^\infty(\Omega)}\leq 1-\delta.
                \end{align*}
                Therefore, thanks to the choice of $\tilde\gamma=\delta$, from \eqref{asympt} and \eqref{ener} we get that
					\begin{align}
						\label{pp}\left({{E}}_{CH}(\vphi(t))-{{E}}_\infty\right)^{1-\vartheta}\leq C\norm{\mu(t)-\overline\mu(t)}_{\Ld} \leq C\| \nabla\mu(t)\|_{\bLd},\quad \forall t\in A_M(T_*),
					\end{align}
                    where of course we can choose $\vartheta<\frac12$ since $\sup_{t\geq 0}E(\varphi(t))<\infty$.
The novel argument introduced in \cite{GPoi} can now be used. Namely, we infer
% core of the novel argument follows. We have, by the energy inequality \eqref{energyineq} and recalling $m(\cdot)\geq m_*>0$ since the mobility is non-degenerate,
% 					\begin{align*}
% 						&\nonumber m_*\int_s^t\norm{\nabla \mu(\tau)}_{\bLd}^2\d\tau \leq \int_s^t \int_\Omega m(\vphi(x,\tau))\normmm{\nabla\mu(x,\tau)}^2\dx\dtau
% 						\leq E(\vphi(s))-E(\vphi(t)),
% 					\end{align*}
% 					for any $t>0$ and almost any $s\in[0,t]$, $s=0$ included.
% 					This gives
% 					\begin{align}
% 						& \left(\int_s^t \norm{\nabla\mu(\tau)}_{\Ld}^2\dtau\right)^{2(1-\vartheta)}
% 						\leq \frac1{m_*^{2(1-\vartheta)}}(E(\vphi(s))-E(\vphi(t)))^{2(1-\vartheta)}.\label{cv}
% 					\end{align}
% 					% Now we observe that
% 					% \begin{align*}
% 						%    \lim_{t\to+\infty} \normmm{\frac12 \norm{\bv(s)}^2_{\bLd}- \frac12 \norm{\bv(t)}^2_{\bLd}}\leq \frac12\normmm{\bv(s)}^2_{\bLd},
% 						% \end{align*}
% 					% by \eqref{veloxconv}.
% 					We now let $t\to\infty$ in \eqref{cv}, and obtain, recalling that $E(\vphi(t))\to E_\infty$ as $t\to \infty$, for almost any $s\in(t_*,\infty)$,
					\begin{align}
						&\nonumber \left(\alpha\int_s^\infty \norm{\nabla\mu(\tau)}_{\Ld}^2\dtau\right)^{2(1-\vartheta)}\\&
						\leq \frac1{m_*^{2(1-\vartheta)}}(E(\vphi(s))-E_\infty)^{2(1-\vartheta)}(\chi_{A_M(T_*)}(s)+\chi_{(T_*,\infty)\setminus A_M(T_*)}(s)),\label{cv22}
					\end{align}
            for almost any $s\in (T_*,\infty)$.
			Observe that, for almost any $s\in (T_*,\infty)\setminus A_M(T_*)$ (i.e., the bad times), it holds $\norm{\nabla \mu(s)}_{\bLd}\geq M$. Therefore, recalling that $E(\vphi(t))\leq E(\vphi_0)$ for any $t\geq0$ and $t\mapsto E(\vphi(t))$ is monotone decreasing, we infer
					$$
				(E(\vphi(s))-E_\infty)^{2(1-\vartheta)}\chi_{(t_*,\infty)\setminus A_M(T_*)}(s)\leq (2E(\vphi_0))^{2(1-\vartheta)}\frac{\norm{\nabla\mu(s)}_{\bLd}^2}{M^2}\chi_{(t_*,\infty)\setminus A_M(T_*)}(s),
				$$
				for almost any $s\in (T_*,\infty)\setminus A_M(T_*)$.
				On the other hand, in the set of good times, thanks to \eqref{pp}, we have
					$$
				(E(\vphi(s))-E_\infty)^{2(1-\vartheta)}\chi_{ A_M(T_*)}(s)\leq C^2{\norm{\nabla\mu(s)}_{\bLd}^2}\chi_{A_M(T_*)}(s),
				$$
				for almost any $s\in  A_M(T_*)$.
				We then deduce from \eqref{cv22} that
					\begin{align}
					& \left(\alpha\int_s^\infty \norm{\nabla\mu(\tau)}_{\bLd}^2\dtau\right)^{2(1-\vartheta)}
					\leq \frac1{m_*^{2(1-\vartheta)}}\left(C^2+\frac{(2 E(\vphi_0))^{2(1-\vartheta)}}{M^2}\right)\norm{\nabla \mu (s)}^2_{\bLd},\label{cv2}
				\end{align}
				for almost any  $s\in (T_*,\infty)$.
			    Using now Lemma \ref{Feireisl} with
$$
Z(\cdot)=\norm{\nabla\mu(\cdot)}_{\bLd}, \quad \tilde\alpha=2(1-\vartheta)\in(1,2), \quad \zeta=\tfrac1{(\alpha m_*)^{2(1-\vartheta)}}(C^2+\tfrac{(2 E(\vphi_0))^{2(1-\vartheta)}}{M^2})>0,
$$ and $\mathcal M=(T_*,\infty)$, we get
					\begin{align}
						\nabla \mu\in L^1(T_*,\infty;\mathbf L^2(\Omega)) . \label{bA}
					\end{align}
					Thus, by comparison, we deduce that $\partial_t\vphi\in L^1(T_*,\infty;\Hu')$. Hence, we have
					$$
					\vphi(t)=\vphi(T_*)+ \int_{T_*}^t\partial_t\vphi(\tau) \: \d\tau \to {\vphi_\infty}\quad \text{ in }\Hu', \text{ as } t \to \infty,$$for some ${\vphi_\infty}\in \Hu'$. This entails that $\vphi(t)$ converges in $\Hu'$ as $t\to\infty$ and we conclude that $\omega(\vphi)$ is a singleton.  In order to show \eqref{equil}, it is enough to recall \eqref{convergenceA}, as now $\omega(\vphi)$ is a singleton. The proof is complete.
\subsection{The conserved Allen--Cahn equation}
\subsubsection{Proof of Lemma \ref{convaaab}}
\label{secconvaaa1}
We consider a sequence $t_n\to \infty$ such that $\varphi(t_n)\rightharpoonup  \widetilde{\varphi}$ weakly in $H^1(\Omega)$, with $\widetilde\vphi\in \omega_1(\varphi)$, and we focus on a nonrelabeled subsequence such that $\vphi(t_n)\to\widetilde\vphi$ strongly in $\Ld$. We set $\varphi_n(t):=\varphi(t+t_n)$ and $\mu_n(t):=\mu(t+t_n)$ and we observe that
\begin{align}
      \label{eqs1a}  & \partial_t\vphi_n+\beta(\mu_n-\overline\mu_n)=0,\quad\text{ a.e. in }\Omega\times(0,+\infty),\\&
        \mu_n=-\gamma\Delta\vphi_n+f'(\vphi_n)\quad\text{ a.e. in }\Omega\times(0,+\infty),\label{eqs2a}
    \end{align}
   together with $\partial_{\bn}\varphi_n=0$ almost everywhere on $\partial\Omega \times [0,+\infty)$.
Using Theorem \ref{weakAC}, we get ${{E}}_1(\vphi(t_n))\leq {{E}}_1(\vphi_0)$ for any $n\in\mathbb{N}$, so that, using standard arguments (see, e.g., \cite{GPoi,GPCAC}), for any $T>0$ we can find a constant $C(T)>0$ independent of $n$ such that
				\begin{align}
				\| \vphi_n\|_{L^\infty(0,T;H^1(\Omega))}+\| \mu_n\|_{L^2(0,T;L^2(\Omega))}+  \| \partial_t\vphi_n\|_{L^2(0,T;L^2(\Omega))}\leq C(T).\label{energy1}
				\end{align}
Also, elliptic regularity gives
\begin{align}
\| \vphi_n\|_{L^2(0,T;H^2(\Omega))}\leq C(T).
    \label{elliptic1}
\end{align}
From estimates \eqref{energy}-\eqref{elliptic}, we deduce that there exists $(\vphi^*,\mu^*)$ such that, for any fixed $T>0$,
				\begin{align}
	&
        \label{estC}\vphi_n\rightharpoonup \vphi^*\quad\text{ in } L^2(0,T;H^2(\Omega))\cap H^1(0,T;\Hu'),\\&
					\vphi_n\overset{\ast}{\rightharpoonup} \vphi^*\quad\text{ in } L^\infty(\Omega\times(0,T)),\label{wstar1}
					\\&
					\vphi_n\to \vphi^*\quad \text{ in } L^2(0,T;H^s(\Omega)),\, \forall\: s\in[0,2)\text{ and a.e. in }\Omega\times(0,T),\label{essenziale1}\\&
					\mu_n\rightharpoonup  \mu^*\quad\text{ in } L^2(0,T;\Ld).
				\end{align}
				% Additionally, by the Aubin-Lions Lemma we also infer that, for any $T>0$,
				% \begin{align}
				% 	\vphi_n\to \vphi^*\quad\text{strongly in }C([0,T];L^2(\Omega)).
				% 	\label{AL}
				% \end{align}
				These convergences are enough to pass to the limit in the equations \eqref{eqs1a}-\eqref{eqs2a}, meaning that {the limit pair } $(\vphi^*,\mu^*)$ satisfies
				\begin{align}
      \label{eqs1a2}  & \partial_t\vphi^*+\beta(\mu^*-\overline\mu^*)=0,\quad\text{ a.e. in }\Omega\times(0,+\infty),\\&
        \mu^*=-\gamma\Delta\vphi^*+f'(\vphi^*)\quad\text{ a.e. in }\Omega\times(0,+\infty),\label{eqs2a12}
    \end{align}
				with initial datum $\vphi^*(0)={\vphi_\infty}$ and boundary condition $\partial_{\bn}\varphi^*=0$ almost everywhere on $\partial\Omega$.  Furthermore, from \eqref{wstar1}-\eqref{essenziale1} we deduce that
               \begin{align}
                 \label{strongc1} & \vphi_n(t)\to \vphi^*(t)\quad\text{ in }H^s(\Omega),\quad\forall s\in(0,2),\\&
                    \vphi_n(t)\to \vphi^*(t) \quad\text{ in }L^p(\Omega),\quad\forall p\geq1,
               \end{align}
               as $n\to\infty$, for almost any $t\geq0$. Therefore, it is immediate to infer
               \begin{align}
               \lim_{n\to\infty}E_1(\vphi_n(t))=E_1(\vphi^*(t)),\quad\text{ for a.a. }t\geq0.
                   \label{enerconverg}
               \end{align}
             Also, by the energy inequality \eqref{energyineq2}, we infer that the energy ${{E}}_1(\vphi(\cdot))$ is nonincreasing, thus there exists ${{E}}_{1,\infty}$ such that
                \begin{align}
				\lim_{t\to+ \infty}{{E}}(\vphi(t))={{E}}_{1,\infty}.
				\label{ene}
                \end{align}
				 Hence, {for almost any $t\geq0$}, we have that (see \eqref{enerconverg})
				$${{E}}_1(\vphi^*(t))={{E}}_{1,\infty}.$$
                Therefore, $E_1(\vphi^*(\cdot))$ is constant in time and equal to $E_{1,\infty}$. Passing then to the limit in the energy inequality, we obtain
				
			\begin{align}\label{zeros1}
				{{E}}_{1,\infty}+\beta\int_s^{t}\int_\Omega \normmm{\mu^*(\tau)-\overline{\mu^*}(\tau)}^2\dx \: \d \tau\leq {{E}}_{1,\infty}\quad \text{ for a.a. } 0\leq s\leq t<\infty,\end{align}
			with $s=0$ included.

    Then, $\mu^*=const$ almost everywhere in $\Omega$. On the other hand, by comparison, it also holds $\partial_t\vphi^*=0$, for almost every $x\in\Omega)$ and $t\geq 0$. Therefore, we deduce
    $$\vphi^*(t)={\vphi_\infty}$$ almost everywhere in $\Omega$, for all $t\geq 0$, so that $\mu_\infty$ is also constant in time. Summing up, ${\vphi_\infty}$ satisfies \eqref{conv1tb} for some constant $\mu_\infty\in \R$, that is, ${\vphi_\infty}\in \mathcal{S}_1$. Also, by \eqref{ene} we infer \eqref{einf}.

Then, the uniform strict separation property of $\omega_1(\vphi)$ as well as the uniform bound in $H^2(\Omega)$ can be obtained following word by word the proof of \cite[Lemma 3.3]{GPoi}, as the stationary problem for the equation Allen--Cahn equation is identical to the one for the Cahn--Hiliard equation without diffusion (see \cite[Definition 3.1]{GPoi}).
In addition, $\vphi\in BC([0,+\infty);H^1(\Omega))$ entails \eqref{convergence1a}.

The main novelty of this lemma is the validity of an asymptotic strict separation for the entire weak solution. We now prove this result by means of De Giorgi's iterations (see also \cite{GPCAC}). First, thanks to \eqref{sepaglobal1a} and \eqref{convergence1a}, we can argue as in \cite[Lemma 3.4]{GPoi} to obtain that, letting
    \begin{align}
    	A_\delta(t) &:= \{x\in\Omega:\; |\varphi(x,t)|\geq 1-{\delta_1}\},\quad t\geq0 ,\label{Adelta}
    \end{align}
    it holds
    \begin{align}
    	\lim_{t\to +\infty }\normmm{A_\delta(t)}\to 0,
    	\label{conerg}
    \end{align}
    where $\delta_1>0$ is given in \eqref{sepaglobal1a}.

 Let us now fix $\delta=\delta_1$. Set then $T>0$ and $\widetilde{\tau}>0$ such that $T-3\widetilde{\tau}\geq 0$.
 We now define, as usual in this kind of argument, the sequence
\begin{align}
k_n=1-\delta-\frac{\delta}{2^n}, \quad \forall n\geq 0,
\label{kn}
\end{align}
where
\begin{align}
1-2\delta< k_n<k_{n+1}<1-\delta,\qquad \forall n\geq 1,\qquad k_n\to 1-\delta\qquad \text{as }n\to \infty,
\label{kn1}
\end{align}
and the sequence of times
\begin{align}
\label{relation}
\begin{cases}
t_{-1}=T-3\widetilde{\tau},\\
t_n=t_{n-1}+\frac{\widetilde{\tau}}{2^n},\qquad n\geq 0,
\end{cases}
\end{align}
satisfying
$$
t_{-1}<t_n<t_{n+1}< T-\widetilde{\tau},\qquad \forall n\geq 0.
$$
Consequently, we introduce a cutoff function $\eta_n\in C^1(\R)$ defined by
\begin{align}
\eta_n(t):=\begin{cases}
0,\quad t\leq t_{n-1},\\
1,\quad t\geq t_{n},
\end{cases}\text{ and }\quad \vert \eta^\prime_n(t)\vert\leq \frac{2^{n+1}}{\widetilde{\tau}}.
\label{cutoff}
\end{align}
We then set
\begin{align}
\varphi_n(x,t):=(\varphi-k_n)^+,
\label{phik0}
\end{align}
and, for any $n\geq 0$, we introduce the interval $I_n=[t_{n-1},T]$ and
$$
A_n(t):=\{x\in \Omega: \varphi(x,t)-k_n\geq 0\},\quad \forall t\in I_n.
$$
Clearly, we have
$$
I_{n+1}\subseteq I_n,\qquad \forall n\geq 0,$$
$$A_{n+1}(t)\subseteq A_n(t),\qquad \forall n\geq 0,\qquad \forall t\in I_{n+1}.
$$
We also define
$$
y_n:=\int_{I_n}\int_{A_n(s)}1\dx\d s,\qquad \forall n\geq0.
$$
 Now, for any $n\geq 0$, we take $v=\varphi_n\eta_n^2$
as test function and integrate over $[t_{n-1},t]$, $t_n\leq t\leq T$. This gives
\begin{align}
&\nonumber\int_{t_{n-1}}^t (\partial_t\varphi)\varphi_n\eta_n^2\dx\d s+\gamma\beta\int_{t_{n-1}}^t\int_{A_n(s)}\nabla\varphi\cdot \nabla\varphi_n \eta_n^2 \dx \d s\\&+\beta\int_{t_{n-1}}^t\int_{\Omega}F^\prime(\varphi)(\varphi_n-\overline{\varphi}_n)\eta_n^2\dx\d s =\beta\int_{t_{n-1}}^t\int_{\Omega}\eta_n^2\theta_0\varphi(\varphi_n-\overline\varphi_n)\dx \d s,
\label{phin}
\end{align}
where we used
\begin{align*}
    \int_{t_{n-1}}^t\int_{\Omega}\eta_n^2(\mu-\overline{\mu})\varphi_n\dx\d s=     \int_{t_{n-1}}^t\int_{\Omega}\eta_n^2\mu(\varphi_n-\overline \varphi_n)\dx\d s.
\end{align*}
In this case, since we only have a weak solution at hand, we cannot exploit the regularity $\overline{F'(\vphi)}\in L^\infty(\tau,+\infty)$, for some $\tau>0$ (cf. \cite{GGPC,GPCAC}). Hence, we use a different argument, based on \cite[Proof of Lemma 3.4]{GPoi}. We observe that
\begin{align*}
    &\int_{t_{n-1}}^t\int_{\Omega}F^\prime(\varphi)(\varphi_n-\overline{\varphi}_n)\eta_n^2\dx\d s\\&
    =\int_{t_{n-1}}^t\int_{\{x:\ \varphi_n(x,s)= 0\}}F^\prime(\varphi)(\varphi_n-\overline{\varphi}_n)\eta_n^2\dx\d s+\int_{t_{n-1}}^t\int_{\{x:\ \varphi_n(x,s)>0\}}F^\prime(\varphi)(\varphi_n-\overline{\varphi}_n)\eta_n^2\dx\d s.
\end{align*}
Then, owing to \eqref{conerg}, there exists $\overline{T}=\overline T(\delta_1)$ such that
$$
\normmm{A_\delta(t)}\leq \frac{\normmm{\Omega}}{4},\quad \quad \forall t\geq \overline T.
$$
This means that, recalling $0\leq \varphi_n\leq 2\delta_1$, with $\delta=\delta_1$,
$$
\overline\varphi_n(t)\leq 2\delta\frac{\normmm{A_n(t)}}{\normmm{\Omega}}\leq 2\delta\frac{\normmm {A_\delta(t)}}{\normmm{\Omega}}\leq  \frac{\delta}4,\quad \forall t\geq \overline T.
$$
As a consequence, we have
\begin{align}
0<k_n+\overline\varphi_n(t) \leq 1-\frac{3\delta}4-\frac{\delta}{2^n},\quad \forall t\geq \overline T,
\label{ext}
\end{align}
so that $F'(k_n+\overline\varphi_n(t))\leq F'(1-\frac{3\delta}4)<+\infty$ for any $t\geq \overline T$, since $\delta=\delta_1$ is fixed. We can now argue exactly as in the proof of \cite[Lemma 3.4]{GPoi}, to infer, assuming $T-3\widetilde \tau>\overline T$,
\begin{align}
&\nonumber\int_{t_{n-1}}^t\int_{\Omega}F^\prime(\varphi)(\varphi_n-\overline{\varphi}_n)\eta_n^2\dx\d s\\&\geq -\int_{t_{n-1}}^t\int_{\{x:\ 0<\varphi(x,s)\leq k_n\}}{F^\prime(k_n)\overline{\varphi}_n}\eta_n^2\dx\d s-\int_{t_{n-1}}^tF'(\overline\varphi_n+k_n)\eta_n^2\int_{\{x:\ \varphi(x,s)\leq k_n\}}\overline{\varphi}_n\dx\d s.
    \label{controls}
\end{align}
We now estimate the right-hand side of \eqref{phin}. Note that, as observed in \cite{P}, since $ \vert\varphi\vert<1\text{ almost everywhere in }\Omega$, for any $t\geq 0$, it holds
\begin{align}
\label{delti}
0\leq \varphi_n\leq 2\delta\quad\text{ a.e. in }\Omega,\quad \forall t\geq 0.
\end{align}
Therefore, on account of the boundedness of $\varphi$, we have
\begin{align*}
    &\int_{t_{n-1}}^t\int_{\Omega}\eta_n^2\theta_0\varphi(\varphi_n-\overline\varphi_n)\dx \d s\leq  \int_{t_{n-1}}^t\int_{A_n(s)}\eta_n^2\theta_0\varphi_n\dx\d s+\normmm{\Omega}\int_{t_{n-1}}^t\eta_n^2\theta_0\overline\varphi_n \d s\\&=2\int_{t_{n-1}}^t\int_{A_n(s)}\eta_n^2\theta_0\varphi_n\dx\d s\leq 4\theta_0\delta y_n.
\end{align*}
On the other hand, observe that
\begin{align}
\int_{t_{n-1}}^t(\partial_t\varphi)\varphi_n\eta_n^2\dx \d s=\frac{1}{2}\Vert\varphi_n(t)\Vert^2_{\Ld}-\int_{t_{n-1}}^t\Vert\varphi_n(s)\Vert^2_{\Ld}\eta_n\partial_t\eta_nds,\label{timede}
\end{align}
as well as
\begin{align}
\nonumber
&\int_{t_{n-1}}^t\Vert\varphi_n(s)\Vert_{\Ld}^2\eta_n\partial_t\eta_n\d s=\int_{t_{n-1}}^t\int_{A_n(s)} \varphi_n^2(s)\eta_n\partial_t\eta_n\dx\d s\\&
\leq
\int_{t_{n-1}}^t\int_{A_n(s)} (2\delta)^2\frac{2^{n+1}}{\widetilde{\tau}}\dx\d s\leq \frac{2^{n+3}\delta^2}{\widetilde{\tau}}y_n.
\label{del}
\end{align}
Therefore, using $\varphi_n\in[0,2\delta]$ again, we get
\begin{align*}
&\frac{1}{2}\Vert\varphi_n(t)\Vert^2_{\Ld}+\gamma\beta\int_{t_{n-1}}^t \eta_n^2\Vert \nabla\varphi_n(s)\Vert^2_{\bLd}\d s\\&\leq\left(4\theta_0\delta\beta+\frac{2^{n+3}\delta^2}{\widetilde{\tau}}\right)y_n+\beta\int_{t_{n-1}}^t\int_{\{x:\ 0<\varphi(x,s)\leq k_n\}}{F^\prime(k_n)\overline{\varphi}_n}\eta_n^2\dx\d s\\&\quad+\beta\int_{t_{n-1}}^tF'(\overline\varphi_n+k_n)\eta_n^2\int_{\{x:\ \varphi(x,s)\leq k_n\}}\overline{\varphi}_n\dx\d s\\&
\leq \left(4\theta_0\delta\beta+\frac{2^{n+3}\delta^2}{\widetilde{\tau}}\right)y_n+\beta\left(F'(k_n)+F'(1-\frac{3\delta}4-\frac{\delta}{2^n})\right)\delta y_n\\&
\leq 2^n\left(4\theta_0\delta\beta+\frac{2^{3}\delta^2}{\widetilde{\tau}}+\delta\beta F'(1-\delta)+\delta \beta F'(1-\frac{3\delta}4)\right)y_n:=2^n K_\delta(\widetilde\tau) y_n,
\end{align*}
where we also used \eqref{ext}. Recall that $\delta=\delta_1$. Thus, we get
\begin{align}
\max_{t\in I_{n+1}}\Vert\varphi_n(t)\Vert^2_{\Ld}\leq 2X_n,\qquad  {\color{black}}\int_{I_{n+1}}\Vert \nabla\varphi_n\Vert^2_{\bLd}ds \leq X_n,
\label{est1}
\end{align}
where
$$
X_n:= 2^{n}K_\delta(\widetilde\tau)\max\left\{1,\frac1{\gamma\beta}\right\}y_n:=2^{n}C_\delta(\widetilde\tau)y_n.
$$
On the other hand, for any $t\in I_{n+1}$ and for almost any $x\in A_{n+1}(t)$, we get
\begin{align}
&\nonumber\varphi_n(x,t)=\varphi(x,t)-\left[1-\delta-\frac{\delta}{2^n}\right]\\&=
\underbrace{\varphi(x,t)-\left[1-\delta-\frac{\delta}{2^{n+1}}\right]}_{\varphi_{n+1}(x,t)\geq 0}+\delta\left[\frac{1}{2^{n}}-\frac{1}{2^{n+1}}\right]\geq \frac{\delta}{2^{n+1}},\label{basic}
\end{align}
which implies
\begin{align*}
\int_{I_{n+1}}\int_{\Omega}\vert\varphi_n\vert^3\dx\d s\geq \int_{I_{n+1}}\int_{A_{n+1}(s)}\vert\varphi_n\vert^3\dx\d s\geq \left(\frac{\delta}{2^{n+1}}\right)^3\int_{I_{n+1}}\int_{A_{n+1}(s)}\dx\d s=\left(\frac{\delta}{2^{n+1}}\right)^3y_{n+1}.
\end{align*}
Then, we have
\begin{align}
&\nonumber\left(\frac{\delta}{2^{n+1}}\right)^3y_{n+1}\leq \int_{I_{n+1}}\int_{\Omega}\vert\varphi_n\vert^3\dx\d s\\&=\int_{I_{n+1}}\int_{A_n(s)}\vert\varphi_n\vert^3\dx\d s \leq \left(\int_{I_{n+1}}\int_{\Omega}\vert\varphi_n\vert^{\frac{10}{3}}\dx\d s\right)^{\frac{9}{10}}\left(\int_{I_{n+1}}\int_{A_n(s)}1\ \dx\d s\right)^{\frac{1}{10}}.
\label{est2}
\end{align}
Therefore, by Gagliardo-Nirenberg's inequalities, we get
\begin{align*}
&\int_{I_{n+1}}\int_{\Omega}\vert\varphi_n\vert^{\frac{10}{3}}\dx\d s\leq \hat{C}\int_{I_{n+1}} \left(\Vert\varphi_n\Vert^2_{\Ld}+\Vert\nabla\varphi_n\Vert^2_{\Ld}\right)\Vert\varphi_n\Vert^{\frac{4}{3}}_{\Ld}\d s,
\end{align*}
so that
\begin{align}
\nonumber&\int_{I_{n+1}}\int_{\Omega}\vert\varphi_n\vert^{\frac{10}{3}}\dx\d s\leq \hat{C}\max_{t\in I_{n+1}}\Vert\varphi_n(t)\Vert_{\Ld}^{\frac{4}{3}}\left(3\widetilde \tau\max_{t\in I_{n+1}}\Vert\varphi_n(t)\Vert^2_{\Ld}+\int_{I_{n+1}}\Vert\nabla\varphi_n\Vert^2_{\Ld}\d s\right)\\&
\leq 2^\frac23{{\color{black}}\hat{C}}X_n^{\frac{2}{3}}\left(6\widetilde\tau X_n+\int_{I_{n+1}}\Vert\nabla\varphi_n\Vert^2_{\bLd}\d s\right)\leq 2^\frac23{{\color{black}}\hat{C}}X_n^{\frac{5}{3}}(6\widetilde\tau+1),\label{F1}
\end{align}
and, by \eqref{est1}, we get
\begin{align*}
\int_{I_{n+1}}\int_{\Omega}\vert\varphi_n\vert^{\frac{10}{3}}\dx\d s\leq
2^\frac23{{\color{black}}\hat{C}(6\widetilde\tau+1)}{2^{\frac 5 3n }}C_\delta(\widetilde\tau)^\frac53 y_n^{\frac{5}{3}}.
\end{align*}
From \eqref{est2}, we thus immediately infer
\begin{align}
&\nonumber\left(\frac{\delta}{2^{n+1}}\right)^3y_{n+1}\leq {2^{{\frac{3}{2}n}+\frac35}\hat{C}^{\frac{9}{10}}(6\widetilde\tau+1)^{\frac{9}{10}}C_\delta(\widetilde\tau)^{\frac32}}y_n^{\frac{8}{5}}.
\end{align}
In conclusion, we end up with
\begin{align}
y_{n+1}\leq \frac{2^{{\frac{9}{2}n}+\frac{18}5}\hat{C}^{\frac{9}{10}}(6\widetilde\tau +1)^{\frac{9}{10}}C_\delta(\widetilde\tau)^{\frac32}}{\delta^3}y_n^{\frac{8}{5}},\qquad \forall n\geq 0.
\label{last0}
\end{align}
Thus we can apply Lemma \ref{conv}. More precisely, setting
$$
b=2^\frac{9}{2}>1, \quad C=\frac{2^{\frac{18}5}\hat{C}^{\frac{9}{10}}(6\widetilde\tau +1)^{\frac{9}{10}}C_\delta(\widetilde\tau)^{\frac32}}{\delta^3}>0, \quad \varepsilon=\frac{3}{5},
$$
we have that ${y}_n\to 0$, as long as
$$
{y}_0\leq C^{-\frac{5}{3}}b^{-\frac{25}{9}},
$$
that is,
\begin{align}
y_0\leq \widetilde C_{\delta}(\widetilde\tau),
\label{last2}
\end{align}
where $ C_{\delta}(\widetilde\tau)>0$ is a suitable constant. Then, from \eqref{conerg}, up to enlarging $\overline T$, we have for any $\xi>0$, since $\delta=\delta_1$,
\begin{equation*}
y_0=\int_{I_0}\int_{A_0(s)}1\dx\d s\leq\int_{I_0}\int_{\{x\in\Omega:\ \varphi(x,t) \geq 1-2\delta\}}1\dx\d s
\leq 3\widetilde\tau\xi.
\end{equation*}
Therefore, if we choose $\xi$ sufficiently small (and thus we fix $\overline T(\xi)$) so that
$$
3\widetilde\tau\xi\leq \widetilde C_{\delta}(\widetilde\tau),
$$
then \eqref{last2} holds, where we set (see relation \eqref{relation}) $T-3\widetilde{\tau}\geq \overline T$.
In the end, taking the limit of $y_n$ as $n\to\infty$, we obtain
$$
\Vert(\varphi-(1-\delta_1))^+\Vert_{L^\infty(\Omega\times({T}-\widetilde{\tau},{T}))}=0,
$$
where, thanks to the choice of $T$, $({T}-\widetilde{\tau},{T})\subset (\overline T,+\infty)$.
 We now repeat the very same argument for the case $(\varphi-(-1+\delta))^-$ (using $\varphi_n(t)=(\varphi(t)+k_n)^-$). As a consequence we showed that there exists $\overline T>0$ such that
 \begin{align}
-1+\delta_1\leq \varphi\leq 1-\delta_1 \quad\text{ a.e. in }\Omega\times ({T}-\widetilde{\tau},{T}),\quad T>\overline T+3\widetilde\tau.
\label{endp}
\end{align}
Finally, we can iterate the same procedure in $({T},T+\widetilde{\tau})$, reaching eventually the whole interval $[ T,+\infty)$, with $T>\overline T+3\widetilde\tau$. The proof of \eqref{sepac} is thus concluded.
\subsubsection{Proof of Lemma \ref{twoparts1}}
\label{sec_prooftwoparts1}
Recall that, for fixed $T>0$, we have
\begin{align}
    \sup_{t\in A_M(T)}\norm{\mu(t)-\overline{\mu}(t)}_{\bLd}\leq M.\label{M1b}
\end{align}
Then, using again an argument based on \cite{Kenmochi} (cf. \eqref{ctr}), and recalling that $\vphi\in[-1,1]$ and $\vphi\in L^\infty(0,+\infty;H^1(\Omega))$, we find
\begin{align*}
          \normmm{\int_\Omega{\mu(t)}\dx}\leq C\int_\Omega\normmm{f'(\vphi(t))}\dx\leq C(1+\norm{\mu(t)-\overline\mu(t)}_{\bLd}),
                \end{align*}
which entails that (see \eqref{M1b})
\begin{align}
    \sup_{t\in A_M(T)}\norm{\mu(t)}_{L^2(\Omega)}\leq C(1+M).\label{M2b}
\end{align}
It is now easy to prove (for instance by testing formally the equation defining $\mu$ with $F'(\vphi)$) that
\begin{align}
    \sup_{t\in A_M(T)}\norm{F'(\vphi(t))}_{L^2(\Omega)}\leq C(1+M),\label{M2bA}
\end{align}
and, by standard elliptic regularity, we infer
\begin{align}
   \sup_{t\in A_M(T)} \norm{\vphi(t)}_{H^2(\Omega)}\leq C(1+M),\quad \forall t\in A_M(T),\label{prb}
\end{align}
and for any $T>0$.

Therefore, the same contradiction argument used in the proof of Lemma \ref{twoparts}, and the embedding $H^2(\Omega)\hookrightarrow\hookrightarrow H^r(\Omega)$, $r\in(0,2)$, yield \eqref{precomp1a} as well as the compactness properties of $\omega_{g,1}(\vphi)$. This ends the proof.

\subsubsection{Proof of Theorem \ref{uniqueeq1a}}
\label{sec:proofLoja1}
The proof can be carried out exactly as the one of Theorem \ref{uniqueeq} in Section \ref{sec:proofLoja}, as long as one replaces $\nabla\mu$ with $\mu-\overline\mu$, and uses the \L ojasiewicz--Simon inequality of Proposition \ref{Lojaw} in place of the one of Proposition \ref{Lojaw2}. This allows the precompactness of trajectories along sequences of good times in $H^1(\Omega)$, given by \eqref{precomp1a}. This suffices to complete the argument, arguing as  in the proof of Theorem \ref{uniqueeq} with $H^1(\Omega)$ in place of $H^2(\Omega)$. This gives (cf. \eqref{bA})
\begin{align}
						 \mu-\overline\mu\in L^1(T_*,+\infty;\mathbf L^2(\Omega)),  \label{bA1}
					\end{align}
                    for some $T_*>0$ sufficiently large.
					Thus, by comparison, we have  $\partial_t\vphi\in L^1(T_*,+\infty;\Ld)$, so that
					$$
					\vphi(t)=\vphi(T_*)+ \int_{T_*}^t\partial_t\vphi(\tau) \: \d\tau \to {\vphi_\infty}\quad \text{ in }\Ld, \text{ as } t \to \infty,$$for some ${\vphi_\infty}\in \Ld$. Therefore, $\omega(\vphi)$ is a singleton. The convergence \eqref{equil1a} is then a straightforward consequence of  \eqref{convergence1a}. The proof is finished.
\subsection{The nonlocal Cahn--Hilliard equation with non-degenerate mobility}
\subsubsection{Proof of Lemma \ref{convaaa3}}
\label{secconvaaa2}
We consider a sequence $t_n\to +\infty$ such that $\varphi(t_n)\rightharpoonup  \widetilde{\varphi}$ weakly in $L^2(\Omega)$, with $\widetilde\vphi\in \omega_2(\varphi)$. Without loss of generality we assume $\vphi(t_n)\to\widetilde\vphi$ strongly in $\Hu'$. We then set $\varphi_n(t):=\varphi(t+t_n)$ and $\mu_n(t):=\mu(t+t_n)$ and we observe that
\begin{align}
      \label{eqs1c}  &\langle \partial_t\vphi_n,v\rangle_{\Hu',\Hu}+\alpha(m(\vphi_n)\nabla \mu_n,\nabla v)=0,\quad \forall v\in H^1(\Omega), \text{ for a.a. }t\geq 0,\\&
        \mu_n=-J\ast\vphi_n+F'(\vphi_n)\quad\text{ a.e. in }\Omega\times(0,+\infty).\label{eqs2c}
    \end{align}
Thanks to Theorem \ref{weaknonloc}, on account of the energy inequality \eqref{energyineq2}, we get that ${{E}}_2(\vphi(t_n))\leq {{E}}_2(\vphi_0)$ for any $n$. Thus, for any $T>0$, there exists $C(T)>0$ independent of $n$ such that
				\begin{align}
				\| \mu_n\|_{L^2(0,T;H^1(\Omega))}\leq C(T).\label{energyc}
				\end{align}
                This term can be controlled, for instance, using the energy inequality and the Poincaré inequality, recalling that (see, e.g., \cite{Kenmochi, FrigeriGrasselli})
                $$
\int_\Omega\normmm{F'(\vphi_n)}\dx\leq C(\overline\vphi_0)\left(1+\int_\Omega F'(\vphi_n)(\vphi_n-\overline\vphi_n)\dx\right),
                $$
                which gives
                \begin{align}
                   \int_\Omega\normmm{\mu_n}\dx=\int_\Omega\normmm{f'(\vphi_n)}\dx\leq C(1+\norm{\nabla\mu_n}_{\bLd}).\label{ctrc}
                \end{align}
By comparison, we also get
\begin{align}
    \| \partial_t\vphi_n\|_{L^2(0,T;H^1(\Omega)')}\leq C(T).\label{cmpc}
\end{align}
Furthermore, testing formally (but it can be made rigorous using suitable approximations) the equation defining $\mu_n$ with $F'(\vphi_n)$ and recalling \eqref{ctrc}, we get
$$
\norm{F'(\vphi_n)}_{L^2(0,T;L^2(\Omega))}\leq C(1+\norm{\nabla\mu}_{L^2(0,T;\bLd)})\leq C(T).
$$
Then, taking (formally) the scalar product of $\nabla\mu_n$   with $\nabla\vphi_n$ and integrating over $\Omega$, we find
\begin{align*}
    \int_\Omega F''(\vphi_n)\normmm{\nabla\vphi_n}^2\dx\leq C(\norm{\nabla J}_{L^1(A)}+\norm{\nabla\mu_n}_{\bLd}^2),
\end{align*}
where $A\subset \R^d$ is a sufficiently large compact set containing $\Omega-\Omega$, entailing from \eqref{energyc} that
\begin{align}
    \norm{\vphi_n}_{L^2(0,T;H^1(\Omega))}\leq C(T),\quad \forall n\in\N,\label{rrr}
\end{align}
for any $T>0$.
From estimates \eqref{energyc}-\eqref{rrr}, we deduce that there exist $\vphi^*$ and $\mu^*$ such that, for any fixed $T>0$,
				\begin{align}
	&
        \label{estA1}\vphi_n\rightharpoonup \vphi^*\quad\text{ in } L^2(0,T;\Hu)\cap H^1(0,T;\Hu'),\\&
					\vphi_n\overset{\ast}{\rightharpoonup} \vphi^*\quad\text{ in } L^\infty(\Omega\times(0,T)),
					\\&
					\vphi_n\to \vphi^*\quad \text{ in } L^2(0,T;H^s(\Omega)),\, \forall\: s\in[0,1)\text{ and a.e. in }\Omega\times(0,T),\label{essenzialek}\\&
					\mu_n\rightharpoonup  \mu^*\quad\text{ in } L^2(0,T;\Hu).
				\end{align}
				Additionally, by \eqref{essenzialek}, we infer
				\begin{align}
					\vphi_n(t)\to \vphi^*(t)\quad\text{strongly in }H^s(\Omega),\quad s\in(0,1),
					\label{ALk}
				\end{align}
                for almost any $t\geq0$.
				These convergences are enough to pass to the limit in the equations \eqref{eqs1c}-\eqref{eqs2c}, proving that the pair $(\vphi^*,\mu^*)$ satisfies, for any $T>0$,
				\begin{align*}
					\langle \partial_t\vphi^*,v\rangle_{H^1(\Omega)',H^1(\Omega)}+\alpha(m(\vphi^*)\nabla \mu^*, \nabla v)=0,\quad \forall \: v\in \Hu,\quad \text{a.e. in }(0,T),\\
					\mu^*=-J\ast \vphi^*+F'(\vphi^*),\quad\text{a.e. in }\Omega\times(0,T),
				\end{align*}
				with initial datum $\vphi^*(0)={\vphi_\infty}$. Furthermore, we have by \eqref{ALk} that
                \begin{align}
				\lim_{n\to \infty}{{E}}_2(\vphi_n(t))={{E}}_2(\vphi^*(t))
				\label{A11}
                \end{align}
                for {almost any} $t\geq 0$. The energy inequality \eqref{energyineq3} implies that ${{E}}_2(\vphi(\cdot))$ is nonincreasing. Thus, there exists ${{E}}_{2,\infty}$ such that
                \begin{align}
				\lim_{t\to+ \infty}{{E}}_2(\vphi(t))={{E}}_{2,\infty}.
				\label{enec}
                \end{align}
				 Hence, {for almost any $t\geq0$}, we have
				$${{E}}_2(\vphi^*(t))=\lim_{n\to \infty}{{E}}(\vphi(t+t_n))={{E}}_{2,\infty}.$$
                Passing then to the limit in the energy inequality (see \eqref{energyineq3}), we obtain
				
			\begin{align}\label{zerosc}
				{{E}}_{2,\infty}+\alpha\int_s^{t}\int_\Omega m(\vphi^*(\tau))\normmm{\nabla\mu^*(\tau)}^2\dx \: \d \tau\leq {{E}}_{2,\infty}\quad \text{ for a.a. } 0\leq s\leq t<+\infty,\end{align}
			with $s=0$ included.

    Then, being $m(\cdot)\geq m_*>0$, \eqref{zerosc} gives $\mu^*=const$ almost everywhere in $\Omega$ and, by comparison, we infer that $\partial_t\vphi^*=0$ in $\Hu'$, for almost every $t\geq 0$. Thus, $$\vphi^*(t)={\vphi_\infty}$$ almost everywhere in $\Omega$, for all $t\geq 0$. This allows us to conclude that $\mu_\infty$ is also independent of time. Therefore, ${\vphi_\infty}$ satisfies \eqref{conv1tb1} for some constant $\mu_\infty\in \R$, and then ${\vphi_\infty}\in \mathcal{S}_2$. Also, \eqref{Einfty} holds due to \eqref{A11}-\eqref{enec}.
Then, precompactness of trajectories \eqref{convergenceAA} can be obtained recalling that $\vphi\in BC([0,+\infty);L^2(\Omega))$ and $L^2(\Omega)\hookrightarrow\hookrightarrow (H^s(\Omega))'$ for any $s\in(0,1]$.
In conclusion, the further characterization \eqref{omegal3} is a direct consequence of \eqref{ALk}.

    % A straightforward consequence of the precompactness of trajectories is then property \eqref{convergence}, since, if this does not hold, by contradiction we would have that there exists $\varepsilon>0$ and a sequence $\{t_n\}_{n\in\N}$, with $t_n\to\infty$, such that
    % \begin{align}
    % 	\inf_{\varphi_\infty\in \omega(\varphi)}\norm{\varphi(t_n)-\varphi_\infty}_{\Hu}>\varepsilon,\quad \forall n\in\N,\label{ess}
    % \end{align}
    % but $\{\varphi(t_n)\}_{n\in\N}$ is uniformly bounded in $H^1(\Omega)$, thus there exists a (nonrelabeled) subsequence such that $\varphi(t_n)\rightharpoonup \varphi_\infty$,  in $H^1(\Omega)$ for some $\varphi_\infty\in H^1(\Omega)$, so that $\vphi_\infty \in \omega(\vphi)$. The previous argument (see \eqref{H1conv3}) entails that there exists a subsequence such that $\vphi(t_n)\to \vphi_\infty$ strongly in $H^1(\Omega)$, which contradicts \eqref{ess}. This also immediately gives that $\omega(\vphi)$ is compact in $H^1(\Omega)$.

% We now show the uniform strict separation properties of the $\omega$-limit.  First, given $\varphi_\infty\in \mathcal S$, there exists $\delta_{\vphi_\infty}>0$, possibly depending on $\varphi_\infty$ such that
% \begin{align*}
%     \norm{\vphi_\infty}_{L^\infty(\Omega)}\leq 1-\delta_{\vphi_\infty},
% \end{align*}
% which is trivially seen from \eqref{conv1t}, since we have $\norm{F'(\vphi_\infty)}_{L^\infty(\Omega)}\leq C(1+\normmm{\mu_\infty})$.
We now show that  $\omega_2(\vphi)$ is bounded in $H^1(\Omega)$. More precisely, arguing again formally (however, see, e.g., \cite{HeWu2}), taking the gradient of \eqref{conv1tb1}, recalling that $\mu_\infty$ is a constant, and testing it by $\nabla{\vphi_\infty}$, we get
					\begin{align*}
					&\theta \|\nabla {\vphi_\infty}\|^2_{\bLd}\leq (\nabla J\ast {\vphi_\infty},\nabla{\vphi_\infty})\\
                    &\leq \| \nabla J\|_{L^1(A)}\| {\vphi_\infty}\|_{L^\infty(\Omega)}\| \nabla{\vphi_\infty}\|_{\bLd}\leq \| \nabla J\|_{L^1(A)}\| \nabla{\vphi_\infty}\|_{\bLd},
					\end{align*}
                    which yields
					\begin{align}
					    \|\nabla {\vphi_\infty}\|\leq  \frac {1} \theta\|\nabla J\|_{L^1(A)} .\label{alig}
					\end{align}
				Therefore, $\omega_2(\vphi)$ is also bounded in $H^1(\Omega)$.
Moreover, the characterization \eqref{omegal3} entails that $\omega_2(\vphi)$ is closed in $\Ld$, and, being bounded in $H^1(\Omega)$, it is compact in $L^2(\Omega)$.

Finally, a simple contradiction argument, recalling that $\vphi\in BC([0,+\infty);L^2(\Omega))$ and using the compact embedding $\Ld\hookrightarrow\hookrightarrow (H^s(\Omega))'$ for any $s\in(0,1]$, allows us to deduce \eqref{convergenceAA}. The proof is concluded.
\subsubsection{Proof of Lemma \ref{twoparts2}}
\label{sec_prooftwoparts2}
Recall that, on the set of good times $A_{M,2}(T)$, for fixed $T>0$, we have
\begin{align}
    \sup_{t\in A_{M,2}(T)}\norm{\nabla\mu(t)}_{\bLd}\leq M.\label{M1b1}
\end{align}
Also, recalling \eqref{ctrc}, we get (see \eqref{M1b1})
\begin{align}
    \sup_{t\in A_{M,2}(T)}\norm{\mu(t)}_{H^1(\Omega)}\leq C(1+M).\label{M2b2}
\end{align}
Arguing as above (see \eqref{alig}), we deduce, by Cauchy-Schwarz inequality,
\begin{align*}
    &\theta\norm{\nabla\vphi(t)}_{\bLd}^2\\&\leq \norm{\nabla J}_{L^1(A)}+\norm{\nabla \mu(t)}_{\bLd}\norm{\nabla\vphi(t)}_{\bLd}\\&\leq \norm{\nabla J}_{L^1(A)}+\frac2\theta\norm{\nabla \mu(t)}_{\bLd}^2+\frac\theta2\norm{\nabla\vphi(t)}_{\bLd}^2,
\end{align*}
where $A\subset \R^d$ is a compact set containing $\Omega-\Omega$.
Thus, using \eqref{M2b2}, we find
\begin{align}
\sup_{t\in A_{M,2}(T)}\norm{\vphi(t)}_{H^1(\Omega)}\leq C(1+M).
    \label{ctrlx}
\end{align}
Therefore, by the same contradiction argument used in the proof of Lemma \ref{twoparts} and exploiting the compact embedding $H^1(\Omega)\hookrightarrow\hookrightarrow H^r(\Omega)$, $r\in(0,1)$, we can infer \eqref{precomp2} as well as the compactness properties of $\omega_{g,2}(\Omega)$ in $H^r(\Omega)$, with $r\in(0,1)$.

We are left to use the notion of good times to prove that the $\omega$-limit set is uniformly strictly separated from the pure phases. Indeed, given $\vphi_\infty\in \omega_{2}(\vphi)$, and $\{t_n\}$ with $t_n\to\infty$, such that $\vphi(t_n)\rightharpoonup \vphi_\infty$ weakly in $\Ld$, arguing as in the proof of Lemma \ref{convaaa3}, we get, up to a nonrelabeled subsequence,
\begin{align*}
    \mu(\cdot+t_n)\rightharpoonup \mu_\infty=\overline{-J\ast \vphi_\infty+f'(\vphi_\infty)} \quad \text{ in }L^q(0,T;H^1(\Omega)),
\end{align*}
for any $q\in[1,2]$, for some constant $\mu_\infty$. Thus, fixing $M>0$, using \eqref{M2b2} and lower semi-continuity, we get
\begin{align}   {T^\frac23\normmm{\Omega}^\frac23}\normmm{\mu_\infty}\leq \liminf_{n\to\infty}\left(\int_{0}^T\norm{\mu(t+t_n)}_{H^1(\Omega)}^\frac32\dt\right)^\frac23\leq  T^\frac23 C(1+M), \label{muc}
\end{align}
where in the last inequality we used the fact that (cf. the definition of $A_{M,2}(t_n)$)
\begin{align}
\nonumber&\int_{0}^T\norm{\mu(t+t_n)}_{H^1(\Omega)}^\frac32\dt=\int_{t_n}^{T+t_n}\norm{\mu(t)}_{H^1(\Omega)}^\frac32\dt\\&\nonumber=\int_{(t_n,T+t_n)\cap A_{M,2}(t_n)}\norm{\mu(t)}_{H^1(\Omega)}^\frac32\dt+\int_{(t_n,T+t_n)\setminus A_{M,2}(t_n)}\norm{\mu(t)}_{H^1(\Omega)}^\frac32\dt\\&\nonumber
    \leq
TC^\frac32(1+M)^\frac32+\normmm{(t_n,T+t_n)\setminus A_{M,2}(t_n)}^\frac14\left(\int_{t_n}^{T+t_n}\norm{\mu(t)}_{H^1(\Omega)}^2\dt\right)^\frac34\\&
\leq TC^\frac32(1+M)^\frac32+ C(T)\normmm{(t_n,T+t_n)\setminus A_{M,2}(t_n)}^\frac14.
\label{from}
\end{align}
Here we used \eqref{energyineq3} and \eqref{ctrc}. Then, using \eqref{energyineq3} again, we get
\begin{align*}
\normmm{(t_n,T+t_n)\setminus A_{M,2}(t_n)}\leq \frac{\int_{t_n}^{T+t_n}\norm{\nabla \mu(t)}^2_{\bLd}\dt}{M^2} \leq \frac{\int_{t_n}^{\infty}\norm{\nabla \mu(t)}^2_{\bLd}\dt}{M^2}\to 0,\text{ as }n\to\infty,
\end{align*}
so that (see \eqref{from})
\begin{align*}
    \liminf_{n\to\infty} \left(\int_{0}^T\norm{\mu(t+t_n)}_{H^1(\Omega)}^\frac32\dt\right)^\frac23\leq T^\frac23(1+M).
\end{align*}
Observing that the right-hand side of \eqref{muc} does not depend on $\vphi_\infty$, and $\normmm{\Omega}$, we get
\begin{align*}
    \sup_{\vphi_\infty\in \omega_{2}(\vphi)}\normmm{\overline{-J\ast \vphi_\infty+f'(\vphi_\infty)}}\leq C(M).
\end{align*}
Therefore, the uniform strict separation property follows by comparison in \eqref{conv1tb1}.
Indeed, one gets
$$ \sup_{\vphi_\infty\in \omega_{2}(\vphi)}\norm{F'(\vphi_\infty)}_{L^\infty}\leq C(M),$$
which gives \eqref{sepaglobal3}, $F'$ being singular at the pure phases. This ends the proof.
\subsubsection{Proof of Theorem \ref{sep}}
\label{sec::proofsep}
To prove validity of the asymptotic strict separation, we first need to show that, for any $\varepsilon>0$, there exists $\tilde T=\tilde T(\varepsilon)>0$ such that
\begin{align}
    \normmm{A_\delta(t)}\leq \varepsilon,\quad \forall t\in A_{M,2}(\tilde T),\label{c1b}
\end{align}
where
\begin{align*}
    A_\delta(t):=\{x\in \Omega: \normmm{\vphi(x,t)}\geq 1-\delta_2\}.
\end{align*}
Now, by the very same argument used to get \cite[Eq. (5.2)]{GPoi}, thanks to \eqref{sepaglobal3}, we get, for any $t\geq0$,
\begin{align*}
    \normmm{A_\delta(t)}\leq \frac{1}{\delta_2^2}\inf_{\vphi_\infty\in \omega_{g,2}(\vphi)}\norm{\vphi(t)-\vphi_\infty}_{\Ld}^2=\frac 1{\delta_2^2}\dist_{\Ld}(\vphi(t),\omega_{g,2}(\vphi))^2,
\end{align*}
which gives \eqref{c1b} by means of \eqref{precomp2}. Notice that, differently from the local Cahn-Hilliard and conserved Allen-Cahn cases, property \eqref{c1b} holds \textit{only }on the set of good times, due to the lack of $L^2(\Omega)$-precompactness of the trajectories for general times. This introduces nontrivial issues in the next argument.

We can now perform a De Giorgi's iteration scheme as in \cite{P}. In particular, we can use the same notation \eqref{kn}-\eqref{phik0} as in the proof of Lemma \ref{convaaab}, choosing $\delta=\delta_2$. For any $n\geq 0$, we consider the test function $v=\vphi_n\eta_n^2$ in \eqref{weakformul}, and integrate over $[t_{n-1},t]$, $t_n\leq t\leq T$. This gives
\begin{align}
&\nonumber\int_{t_{n-1}}^t\langle\partial_t\vphi,\vphi_n\eta_n^2\rangle_{H^1(\Omega)',\Hu}ds+\alpha\int_{t_{n-1}}^t\int_{A_n(s)}m(\vphi)F^{\prime\prime}(\vphi)\nabla\vphi\cdot \nabla\vphi_n \eta_n^2 \dx \d s\\&=\int_{t_{n-1}}^t\int_{A_n(s)}\eta_n^2(\nabla J\ast \vphi)\cdot \nabla \vphi_n\dx \d s.
\label{phin1}
\end{align}
Here we used the fact that $\nabla F^\prime(\vphi(t))=F''(\vphi)\nabla\vphi(t)$, for almost every $x\in\Omega$ and for almost any $t\geq 0$, which can be proven, for instance, by a truncation argument (see, e.g., \cite[Lemma 3.2]{Wu}).
Observe now that
\begin{align}
\alpha\int_{t_{n-1}}^t\eta_n^2\int_{A_n(s)}m(\vphi)F^{\prime\prime}(\vphi)\nabla\vphi\cdot \nabla\vphi_n \dx \d s\geq \alpha\theta m_*\int_{t_{n-1}}^t\eta_n^2\Vert \nabla\vphi_n\Vert^2_{\bLd}ds.
\label{Fss}
\end{align}
The right-hand side of \eqref{phin}, recalling that $\vert\vphi\vert< 1$ almost everywhere in $\Omega\times (0,+\infty)$, can be estimated as follows
\begin{align}
&\nonumber\int_{t_{n-1}}^t\int_{A_n(s)}(\nabla J\ast \vphi)\cdot \nabla \vphi_n\eta_n^2\dx \d s\\&\nonumber\leq \frac{\alpha\theta}{2}\int_{t_{n-1}}^t\eta_n^2\Vert \nabla\vphi_n\Vert^2_{\bLd}ds+\frac{1}{2\alpha\theta}\int_{t_{n-1}}^t\int_{A_n(s)}\eta_n^2\vert \nabla J\ast \vphi\vert^2 \dx\d s\\&\nonumber
\leq
\frac{\alpha\theta}{2}\int_{t_{n-1}}^t\eta_n^2\Vert \nabla\vphi_n\Vert_{\bLd}^2ds+\frac{1}{2\alpha\theta}\int_{t_{n-1}}^t\Vert\nabla J\ast \vphi\Vert_{L^\infty(\Omega)}^2\int_{A_n(s)}1\dx\d s\\&\nonumber\leq
\frac{\alpha\theta}{2}\int_{t_{n-1}}^t\eta_n^2\Vert \nabla\vphi_n\Vert^2ds+\frac{\Vert\nabla J\Vert_{{L^1(A)}}^2}{2\alpha\theta}\int_{t_{n-1}}^t\int_{A_n(s)}1\dx\d s\\&\leq \frac{\alpha\theta}{2}\int_{t_{n-1}}^t \eta_n^2\Vert \nabla\vphi_n\Vert_{\bLd}^2ds+\frac{\Vert\nabla J\Vert_{{L^1(A)}}^2}{2\alpha\theta}y_n,
\label{J01}
\end{align}
where we have applied the well known inequality
\begin{align}
\Vert\nabla J\ast \vphi\Vert_{L^\infty(\Omega)}\leq \Vert\nabla J\Vert_{{L^1(A)}}\Vert\vphi\Vert_{L^\infty(\Omega)}\leq \Vert\nabla J\Vert_{{L^1(A)}}.
\label{J}
\end{align}
Here $A\subset \R^d$ is again a compact set containing $\Omega-\Omega$.
Furthermore, we have
\begin{align}
\int_{t_{n-1}}^t\langle\partial_t\vphi,\vphi_n\eta_n^2\rangle_{\Hu',\Hu}\d s=\frac{1}{2}\Vert\vphi_n(t)\Vert^2_{\Ld}-\int_{t_{n-1}}^t\Vert\vphi_n(s)\Vert^2_{\Ld}\eta_n\partial_t\eta_n\d s.
\label{eq1z}
\end{align}
Then, we have (see \eqref{del})
\begin{align}
&\int_{t_{n-1}}^t\Vert\vphi_n(s)\Vert^2_{\Ld}\eta_n\partial_t\eta_n\d s\leq \frac{2^{n+3}\delta^2}{\widetilde{\tau}}y_n.
\label{del1}
\end{align}
{Plugging \eqref{Fss}, \eqref{J01}, \eqref{eq1z} and \eqref{del} into \eqref{phin}, we find}
\begin{equation*}
\frac{1}{2}\Vert\vphi_n(t)\Vert^2+\frac{\alpha\theta}{2}\int_{t_{n-1}}^t \eta_n^2\Vert \nabla\vphi_n(s)\Vert^2\d s
\leq2^{n}\max\left\{\frac{\Vert\nabla J\Vert_{{L^1(A)}}^2}{2\alpha\theta},\frac{8\delta^2}{\widetilde{\tau}}\right\}y_n:=2^{n}\widetilde{K}_\delta(\widetilde\tau) y_n,
\end{equation*}
for any $t\in[t_n,T]$. Therefore, we get
\begin{align}
\max_{t\in I_{n+1}}\Vert\vphi_n(t)\Vert^2\leq X_n,\qquad  \int_{I_{n+1}}\Vert \nabla\vphi_n\Vert^2\d s \leq X_n,
\label{estA}
\end{align}
where
$$
X_n:= 2^{n+1}\widetilde{K}_\delta(\widetilde\tau)\max\left\{\frac1{\alpha\theta},1\right\}y_n:=2^{n+1}\widetilde{C}_\delta(\widetilde\tau)y_n.
$$
We can now repeat estimates \eqref{basic}-\eqref{F1} to obtain
\begin{align}
\int_{I_{n+1}}\int_{\Omega}\vert\varphi_n\vert^{\frac{10}{3}}\dx\d s\leq {{\color{black}}\hat{C}}X_n^{\frac{5}{3}}(3\widetilde\tau +1)y_n^{\frac{5}{3}}\leq 2^{\frac53 n+\frac53}{\widehat C}(3\widetilde \tau+1)\widetilde{C}_\delta(\widetilde\tau)^\frac53y_n^\frac53,\label{F2}
\end{align}
as well as (see \eqref{est2})
\begin{align}
&\nonumber\left(\frac{\delta}{2^{n+1}}\right)^3y_{n+1}\leq 2^{\frac32n+\frac32}{\widehat C}^\frac9{10}(3\widetilde{\tau}+1)^\frac{9}{10}\widetilde{C}_\delta(\widetilde\tau)^\frac32y_n^\frac 85:=2^{\frac{3}2n}\widetilde{C}_{1,\delta}y_n^\frac85(\widetilde\tau).
\end{align}
Summing up, we get
\begin{align}
&\nonumber y_{n+1}\leq \frac{2^{\frac{9}2n+3}}{\delta^3}\widetilde{C}_{1,\delta}y_n^\frac85(\widetilde\tau).
\end{align}
Thus we can apply Lemma \ref{conv} with $b=2^\frac{9}{2}>1$, $C=\frac{2^{3}}{\delta^3}\widetilde{C}_{1,\delta}(\widetilde\tau )>0$, $\varepsilon=\frac{3}{5}$, to get that ${y}_n\to 0$ as $n$ goes to $\infty$, as long as
\begin{align}
y_0\leq \dfrac{2^{-\frac{35}2}\delta^5}{\widetilde{C}_{1,\delta}(\widetilde\tau )^\frac53}.
\label{last}
\end{align}
We are left to prove the last estimate, by making use of \eqref{c1b}, together with the crucial help of the definition of good times. For a fixed $\tilde T >0$, to be chosen later on, setting $I_0=[T-3\widetilde \tau, T]$ so that $I_0\subset [\tilde T,+\infty)$, we have
\begin{align*}
&y_0=\int_{I_0}\int_{A_0(s)}1dx\dt\leq\int_{I_0}\int_{\{x\in\Omega:\ \vphi(x,t) \geq 1-2\delta\}}1dx\dt\\
&\leq \int_{I_0\cap A_M(\tilde T)}\normmm{A_\delta(t)}\dt+\int_{I_0\setminus  A_M(\tilde T)}\normmm{A_\delta(t)}\dt.
\end{align*}
Now, recalling \eqref{c1b}, we choose $\tilde T=\tilde T(\varepsilon)$ so that
\begin{align*}
    \int_{I_0\cap A_M(\tilde T)}\normmm{A_\delta(t)}\dt\leq 3\widetilde\tau\varepsilon,
\end{align*}
whereas, exploiting the definition of the set of good times $A_M(\tilde T)$, we get
\begin{align*}
    &\int_{I_0\setminus  A_M(\tilde T)}\normmm{A_\delta(t)}\dt\leq \normmm{\Omega}\normmm{I_0\setminus A_M(\tilde T)}\\
    &\leq \normmm{\Omega}\int_{I_0\setminus A_M(\tilde T)}\frac{\norm{\nabla \mu(t)}_{\bLd}^2}{M^2}\dt \leq \frac{\normmm{\Omega}}{M^2}\int_{\tilde T}^\infty\norm{\nabla \mu(t)}^2_{\bLd}\dt,
\end{align*}
and, being $\nabla\mu\in L^\infty(0,+\infty;\bLd)$, given $\varepsilon>0$ there exists $\overline T=\overline T(\varepsilon)$ such that
\begin{align*}
    \int_{\overline T}^\infty\norm{\nabla \mu(t)}^2_{\bLd}\dt\leq \varepsilon.
\end{align*}
As a consequence, for every $\varepsilon>0$, we can choose $\widehat T=\max\{\tilde T, \overline T\}$ such that
\begin{align*}
&y_0\leq \left(3\widetilde \tau+\frac{\normmm{\Omega}}{M^2}\right)\varepsilon,
\end{align*}
as long as $I_0\subset [\widehat T,+\infty)$, that is (see \eqref{relation}), if $T-3\widetilde \tau>\widehat  T$.
Therefore, if we now  ensure that $\varepsilon>0$ is so small (and thus $\widehat T=\widehat T(\varepsilon)$ is sufficiently large) such that
$$
 \left(3\widetilde \tau+\frac{\normmm{\Omega}}{M^2}\right)\varepsilon\leq \dfrac{2^{-\frac{35}2}\delta^5}{\widetilde{C}_{1,\delta}(\widetilde\tau )^\frac53}
$$
then \eqref{last} holds. Having fixed $T-3\widetilde \tau>\widehat T$, passing to the limit in $y_n$ as $n\to\infty$, we have thus obtained that
$$
\Vert(\vphi-(1-\delta))^+\Vert_{L^\infty(\Omega\times({T}-\widetilde{\tau},{T}))}=0.
$$
 We can now repeat the same argument for the case $(\vphi-(-1+\delta))^-$ (using $\vphi_n(t)=(\vphi(t)+k_n)^-$),  to get in the end that there exists $\delta=\delta_2>0$ such that
 \begin{align}
-1+\delta\leq \vphi \leq 1-\delta \quad\text{ a.e. in }\Omega\times ({T}-\widetilde{\tau},{T}),
\label{end}
\end{align}
for fixed $T>\widehat T+3\widetilde \tau$.
Finally, since $\delta=\delta_2$ does not depend on the choice of the time interval, we can iterate the same procedure in $({T},T+\widetilde{\tau})$, reaching eventually $[T,+\infty)$. The proof of \eqref{sepac1} is thus concluded.
\subsubsection{Proof of Theorem \ref{uniqueeq2}}
\label{sec:proofLoja2}
The proof is now very similar to the one of Lemma \ref{uniqueeq}, so that we only give the main details. Let us choose $\tilde \gamma$ coinciding with the value of $\delta_B$ (see Theorem \ref{sep}), so that, as $\tilde \gamma\leq\delta_2$, it holds (see \eqref{sepaglobal3}) $-1+\tilde\gamma\leq \vphi_\infty\leq 1-\tilde\gamma$ in $\Omega$, for any $\vphi_\infty\in\omega_{g,2}(\vphi)$. Furthermore, for any $\vphi_{\infty,m}\in \omega_{g,2}(\vphi)$ we can find $\vartheta_m\in \left(0,\frac{1}{2}\right]$ and $\eta_m>0$, given {by }Proposition \ref{Lojaw2}, for which \eqref{Lojaw2} is valid with constant $C_m$. Now, from Lemma \ref{twoparts2} we get that $\omega_{g,2}(\vphi)\subset \omega_2(\vphi)$ is compact in $L^2(\Omega)$, so that we can find a finite family of $L^2(\Omega)$-open balls, say $\{B_{\eta_m}\}_{m=1}^{M_1}$, centered at $\{\vphi_{\infty,m}\}_{m=1}^{M_1}\subset \omega_{g,2}(\vphi)$ and with radii $\eta_m$ (depending on the center $\vphi_{m,\infty}\in\omega_{g,2}(\vphi)$), such that
				\begin{equation*}
					\bigcup_{ \varphi_\infty \in \omega_{g,2}(\varphi)} \{\varphi_\infty\} \subset W:= \bigcup_{m=1}^{M_{1}}B_{\eta_m}.
				\end{equation*}
	Recalling \eqref{E2}, the energy functional $E_2(\cdot)$ is constant over $\omega_2(\vphi)$. Also, as {the centers  $\{\vphi_m\}_{m = 1}^{M_1}$ are in finite number, }we can infer that \eqref{ener1} holds \textit{uniformly}, for any $\vphi\in W$ such that $\|\vphi\|_{L^\infty(\Omega)}\leq 1-\tilde\gamma$, and we can replace $E_2(\vphi_\infty)$ with $E_{2,\infty}$ (where we can choose $\vartheta\in(0,\frac12)$, since $\sup_{t\geq 0}E_2(\varphi(t))<\infty$).

    By \eqref{precomp2}, there exists ${t_*}>0$ such that $\vphi(t)\in W$ for any good time $t\in A_M(t_*)$ and, by \eqref{sepac1}, the uniform strict separation property holds for any $t\geq t_*$, i.e.,
                \begin{align*}
                    \sup_{t\geq t*}\norm{\vphi(t)}_{L^\infty(\Omega)}\leq 1-\delta_B.
                \end{align*}
                As a consequence, thanks to the choice of $\tilde \gamma=\delta_B$, arguing as to get \eqref{cv2}, we obtain
					\begin{align}
					& \left(\alpha\int_s^\infty \norm{\nabla\mu(\tau)}_{\bLd}^2\dtau\right)^{2(1-\vartheta)}
					\leq C(M,E_2(\vphi_0))\norm{\nabla \mu (s)}^2_{\bLd},\label{cv2bA}
				\end{align}
				for almost any  $s\in (t_*,+\infty)$, for a suitable constant $C(M,E_2(\vphi_0))>0$.
			    Using now Lemma \ref{Feireisl} with
$$
Z(\cdot)=\norm{\nabla\mu(\cdot)}_{\bLd}, \quad \tilde\alpha=2(1-\vartheta)\in(1,2), \quad \zeta=\frac{1}{\alpha^{2(1-\vartheta)}}C(M,E_2(\vphi_0))>0,
$$ and $\mathcal M=(t_*,+\infty)$, we find
					\begin{align}
						\nabla \mu\in L^1(t_*,+\infty;\mathbf L^2(\Omega)) . \label{bAA}
					\end{align}
					Thus, by comparison, we deduce that $\partial_t\vphi\in L^1(t_*,+\infty;\Hu')$, so that
					$$
					\vphi(t)=\vphi(t_*)+ \int_{t_*}^t\partial_t\vphi(\tau) \: \d\tau \to {\vphi_\infty}\quad \text{ in }\Hu', \text{ as } t \to +\infty,$$for some ${\vphi_\infty}\in \Hu'$. This means that $\vphi(t)$ converges in $\Hu'$ as $t\to +\infty$. Therefore $\omega(\vphi)$ is a singleton.  In order to show \eqref{equil2b}, it is enough to use \eqref{convergenceAA}. The proof is complete.

                    \section{Conclusions and future developments}
\noindent
Here we have shown how the notion of good times can be extended to discuss good equilibrium points. This allows us to prove the convergence of any global weak solution to a single equilibrium for certain Cahn--Hilliard and Allen--Cahn type equations that are characterized by a precompactness of trajectories which is too weak to apply the \L ojasiewicz--Simon inequality. This completes the first results of \cite{GPoi}, introducing the novel idea that convergence to a single equilibrium, as well as the asymptotic strict separation property (in case of second order PDEs) can be obtained \textit{without} the necessity of  regularization of weak solutions. The present approach can be employed, for instance, to establish the convergence to equilibrium of global weak solutions to some diffuse interface models for binary fluids where the capillary force is not regular enough to split the energy inequality (cf. \cite{GPoi} for the AGG system). Typical examples are nonlocal Cahn--Hilliard--Navier--Stokes systems like the one analyzed in \cite{Frigeri} and the Allen--Cahn--Navier--Stokes systems (see \cite{GGPC,GHP}). Other possible applications are concerned with multi-component Allen--Cahn or Cahn--Hilliard equations (see \cite{GGPS,GPCAC}, see also \cite{AGGmulti}), as well as phase-field models on evolving surfaces (cf. \cite{CEGP,DE}) or subject to dynamic boundary conditions (see, for instance, \cite{GGW,KLLM,LW,Wu}).

	\appendix
	\section{Some technical lemmas}
	\subsection{A Lemma on the integrability of functions}
	The following lemma, whose proof can be found in \cite[Lemma 7.1]{FS}, guarantees that a square integrable function in $(0,+\infty)$ is integrable if it satisfies a suitable integral inequality.		
			
\begin{lemma}
	Let $Z\geq 0$ be a measurable function on $(0,\infty)$ such that
	\begin{align*}
		Z\in L^2(0,+\infty),\quad \int_0^\infty\normmm{Z(t)}^2\dt\leq Y,
	\end{align*}
	for some $Y>0$. If there exist $\tilde\alpha\in (1, 2)$, $\zeta>0$, and an open set $\mathcal M\subset (0,+\infty)$ such that
	\begin{align*}
		\left(\int_s^\infty Z^2(t)\dt\right)^{\tilde\alpha}\leq \zeta Z^2(s),\quad \text{ for a.a. }s\in \mathcal M,
	\end{align*}
	then $Z\in L^1(\mathcal M)$, and there exists $C=C(Y,\tilde\alpha,\zeta)>0$, independent of $\mathcal M$, such that
	\begin{align*}
		\int_{\mathcal M} Z(t)\dt \leq C.
	\end{align*}
	\label{Feireisl}
\end{lemma}
\subsection{A lemma on geometric convergence of sequences}
 This lemma, which is a key tool in De Giorgi's iteration argument, can be found, e.g., in \cite[Ch. I, Lemma 4.1]{DiBenedetto}.
\begin{lemma}
	\label{conv}
	Let $\{y_n\}_{n\in\N\cup \{0\}}\subset \R^+$ satisfy the recursive inequalities
	\begin{align}
		y_{n+1}\leq Cb^ny_n^{1+\varepsilon},
		\label{ineq}\qquad \forall n\geq 0,
	\end{align}
	for some $C>0$, $b>1$ and $\varepsilon>0$. If
	\begin{align}
		\label{condition}
		y_0\leq \theta:= C^{-\frac{1}{\varepsilon}}b^{-\frac{1}{\varepsilon^2}},
	\end{align}
	then
	\begin{align}
		y_n\leq \theta b^{-\frac{n}{\varepsilon}},\qquad \forall n\geq 0,
		\label{yn}
	\end{align}
	and consequently $y_n\to 0$ as $n\to \infty$.
\end{lemma}

%\section{Conclusions and future developments}
\noindent

\medskip
\textbf{Acknowledgments.}
This research was funded in part by the Austrian Science Fund (FWF) \href{https://doi.org/10.55776/ESP552}{10.55776/ESP552}.
AP and MG are also members of Gruppo Nazionale per l’Analisi Ma\-te\-ma\-ti\-ca, la Probabilità e le loro Applicazioni (GNAMPA) of
Istituto Nazionale per l’Alta Matematica (INdAM). This research is part of the activities of “Dipartimento di Eccellenza 2023-2027” of
Politecnico di Milano (MG).
For open access purposes, the authors have applied a CC BY public copyright license to
any author accepted manuscript version arising from this submission.
\\
	%\bibliographystyle{siam}

	%\bibliography{GPoi_weak}
\textbf{Data Availability Statement.} Data sharing not applicable to this article as no
datasets were generated or analyzed during the current study.
\\

\textbf{Declarations
Conflict of interest.} The authors have no relevant financial or non-financial interests
to disclose.

\end{document}